\documentclass[11pt,a4paper]{article}
\usepackage[english]{babel}
\usepackage[latin1]{inputenc}
\usepackage{amssymb}
\usepackage{amsmath}
\usepackage{amsthm}

\newcommand\RR{\mathbb{R}}
\newcommand\tr{\ensuremath{\triangle}}
\newcommand\norm[2]{\ensuremath{|\!| #1 |\!|_{#2}}}
\newcommand\nL[2]{\norm{ #1}{L^{#2}}}
\newcommand\nH[2]{\norm{#1}{H^{#2}}}
\newcommand\nX[2]{\norm{#1}{#2}}
\newcommand\nLx[2]{\norm{#1}{L^{#2}_x}}
\newcommand\nHx[2]{\norm{#1}{H^{#2}_x}}
\newcommand\nHM[2]{\norm{#1}{H^{#2}(M)}}
\newcommand\nLi[3]{\norm{#1}{L^#2(#3)}}
\newcommand\nLd[1]{\norm{#1}{L^2}}
\newcommand\LiL[3]{\ensuremath{L^{#1}({#2},L^{#3}(\mathbb R^d))}}
\newcommand\nLiL[4]{\ensuremath{\norm{#1}{\LiL{#2}{#3}{#4}}}}
\newcommand\nLiLM[4]{\ensuremath{\norm{#1}{L^{#2}(#3,L^{#4}(M))}}}
\newcommand\LL[2]{\ensuremath{L^{#1}(L^{#2})}}
\newcommand\LtLx[2]{\ensuremath{L^{#1}_t L^{#2}_x}}
\newcommand\nLtLx[3]{\ensuremath{\norm{#1}{\LtLx{#2}{#3}}}}
\newcommand\nLL[3]{\ensuremath{\norm{#1}{\LL{#2}{#3}}}}

\newtheorem{thm}{Theorem}[section]

\newtheorem*{citethm}{Theorem}
\newtheorem*{citelem}{Lemma}
\newtheorem{lem}[thm]{Lemma}
\newtheorem{prop}[thm]{Proposition}
\newtheorem{cor}[thm]{Corollary}
\newtheorem*{rem}{Remark}
\newtheorem*{Schur}{Schur's lemma}
\newtheorem{Def}{Definition}

\title{Strichartz Inequalities for Lipschitz Metrics on Manifolds and Nonlinear Schrödinger Equation on Domains}
\author{Ramona Anton}
\date{}
\begin{document}

\maketitle

\begin{abstract}

We prove wellposedness of the Cauchy problem for the nonlinear Schrödinger equation for any defocusing power nonlinearity on a domain of the plane with Dirichlet boundary conditions. The main argument is based on a generalized Strichartz inequality on manifolds with Lipschitz metric.

\end{abstract}

\section{Introduction}
Let $\Omega$ be a compact regular domain of $\RR^d$, where $d=2,3$. 
The problem we are interested in is the Dirichlet problem for the semilinear Schrödinger equation
\begin{equation}
\label{NLS}
\left \{
\begin{array}{rcl}
i \partial_t u + \triangle u &=& |u|^\beta u, \ on \ \mathbb{R} \times
\Omega \\ u_{|_{t=0}} &=& u_0, \ on \ \Omega\\
u_{|_{\mathbb R \times \partial \Omega}} &=& 0.
\end{array}
\right.
\end{equation}
More precisely we are interested in proving global existence results in the energy space $H_0^1(\Omega)$ and this will be done for $d=2$. 

This problem has been extensively study in the case of $\Omega=\RR^d$. Note that the sign of the nonlinearity gives an a priori bound of the $H^1$ norm of the flow and thus allows to prove existence of weak solutions in $C(\RR, H_w^1(\RR^d))$. The existence of global strong solution is more difficult. One of the main ingredient to address this difficulty is the Strichartz inequality for the linear flow $e^{it\tr}$. It can be seen as an improvement of the Sobolev imbedding $H^1 \hookrightarrow L^q$ and the price to pay is an average in time rather than a pointwise information. In $\RR^d$, the Strichartz inequality reads as follows : for $(p,q)$ an admissible pair in dimension $d$ and $u_0\in L^2$
$$\nLiL{e^{it\tr}u_0}{p}{\RR}{q}\leq c\nL{u_0}{2}.$$

Let us recall the definition of an admissible pair.

\begin{Def} A pair $(p,q)$ is called admissible in dimension $d$ if $p\geq 2$, $(p,q,d)\neq (2,\infty,2)$ and $$\frac 2p +\frac dq=\frac d2.$$
\end{Def}

In 1977 Strichartz \cite{STR} proves the particular case $p=q$,
$$\nLi{e^{i\cdot\tr}u_0}{{2+\frac 4 d}}{\RR\times \RR^d}\leq c\nL{u_0}{2}.$$
This was generalized by Ginibre and Velo \cite{GiVe} in 1985 for $L_t^p L_x^q$ norm with $p$ and $q$ that satisfy the admissibility condition with $p>2$ and by Keel and Tao \cite{KT} in 1998 for the endpoint $p=2$ and $q=\frac{2d}{d-2}$. Extension to non homogeneous equation is due to Yajima \cite{Ya} in 1987 and Cazenave and Weissler \cite{CaW} in 1988 : for $(p_1,q_1)$ and $(p_2,q_2)$ admissible pairs and $f\in L^{p_2'}([0,T],L^{q_2'}(\RR^d))$ the solution of the non homogeneous equation $$i \partial_t u + \triangle u =f,\ u_{|_{t=0}} = u_0$$ belongs to $C([0,T],L^2)\cap L^{p_1}([0,T],L^{q_1}(\RR^d))$
and satisfies to
$$\nLiL{u}{p_1}{[0,T]}{q_1} \leq c \nLiL{f}{p_2'}{[0,T]}{q_2'}.$$
A contraction mapping argument and those Strichartz inequalities imply the global existence

\begin{citethm} \textbf{\upshape{(\cite{Ya}, \cite{CaW})}} For $d\geq 2$ and $1\leq\beta <\frac{4}{d-2}$ there exists a unique solution $$u\in C(\RR,H^1(\RR^d))\cap L_{loc}^{p_1}(\RR,W^{1,q_1}(\RR^d)),$$ for each $(p,q)$ admissible pair in dimension $d$, of the equation
$$\left \{
\begin{array}{rcl}
i \partial_t u + \triangle u &=& |u|^\beta u,\\ 
u_{|_{t=0}} &=& u_0.
\end{array}
\right.$$
\end{citethm}

For $\Omega\neq \RR^d$ much less is known. In the case of the tori $\mathbb{T}^d$, $d=2,3$, Bourgain \cite{Bo} proved global existence result using less stringent dispersive estimates. In the case of a boundaryless compact manifold Burq, Gérard and Tzvetkov \cite{BGT1} proved Strichartz inequalities with loss of derivatives and showed that those losses are in some specific geometries unavoidable.

In the case of domains of $\RR^2$ and for cubic equations previous results were proved by Brezis and Gallouet \cite{BrGa} in 1980 and Vladimirov \cite{Vl} in 1984.

\begin{citethm} \textbf{\upshape{(\cite{BrGa}, \cite{Vl})}} For $u_0\in H_0^1(\Omega)$ there exists a unique solution $u\in C(\RR, H_0^1(\Omega))$ of the cubic nonlinear equation 
$$i \partial_t u + \triangle u = |u|^2 u, \ on \ \mathbb{R} \times
\Omega,\ u_{|_{t=0}} = u_0, \ on \ \Omega.$$ Moreover, if $u_0\in H_0^1(\Omega)\cap H^2(\Omega)$ then $$u\in C(\mathbb{R}, H_0^1(\Omega)\cap H^2(\Omega))\cap C^1(\mathbb{R},L^2(\Omega)).$$
\end{citethm}

The main ingredients of the proof are the following logarithmic inequalities.

(\textbf{B.-G.}) $\nX{u}{L^\infty} \leq C\nX{u}{H^1}\left( 1+\log \left( 2+\frac{\nX{u}{H^2}} {\nX{u}{H^1}} \right) \right)^{\frac 12}$.

(\textbf{V.}) $\forall p<\infty$, $\nX{u}{L^p}\leq c\sqrt{p}\nX{u}{H^1}$.

The methods used in this proof do not give us informations about nonlinearities stronger than cubic. Note that even in this cubic case, the proof did not yield the Lipschitz continuity on the energy space, which is a consequence of Strichartz estimate in the case of $\Omega=\RR^d$.

In this article we prove a generalized Strichartz inequality for the Schrödinger flow $e^{it\tr}$, where $\tr$ is the Laplace operator on domains of $\RR^d$, $d=2,3$. Let us introduce the following notation : for every $s\in [0,1]$, we denote by $H_D^s(\Omega)$ the domain of the operator $(-\tr_D)^{\frac{s}{2}}$ in $L^2(\Omega)$, where $\tr_D$ is the Dirichlet Laplacian. We refer to section \ref{dbl} for more details. We translate the problem on the domain into a problem on a boundaryless Riemannian manifold by doing a mirror reflection of the domain and identifying the points on the boundary. We make also an even reflection of the coefficients of the metric over the boundary in normal coordinates. Thus we obtain a metric with Lipschitz coefficients. 

We combine ideas from \cite{BaCh} (see also \cite{Ta}) on regularizing the metric with a semiclassical analysis of the flow like in \cite{BGT1} and obtain the following Strichartz inequality (with loss of derivatives) in a general context: $M$ is a compact (or flat outside a compact set) Riemannian manifold of dimension $d=2,3$, endowed with a Lipschitz metric $G$.

\begin{thm}
\label{StrThm}
Let $I$ be a finite time interval, $(p,q)$ an admissible pair in dimension $d=2,3$. Let $\epsilon>0$ be an arbitrarily small constant. Then there exists a constant $c(p,I)>0$ such that, for all $v_0 \in H^{\frac{3}{2p}+\epsilon}(M)$, the following holds
\begin{equation}
\label{Strest}
\nX{e^{it\tr_G}v_0}{L^p(I,L^q(M))} \leq c(p,I) \nH{v_0}{\frac{3}{2p}+\epsilon}.
\end{equation}
\end{thm}
For a compact $C^2$ perturbation of the Laplacian with nontrapping condition, G.Staffilani and D.Tataru \cite{StTa} proved Strichartz inequalities without loss of derivatives. In $1D$ with $BV$ metric similar results were obtain by V.Banica \cite{VBan}, D.Salort \cite{DSal} and N.Burq and F.Planchon \cite{BuPl}. C.Carlos and E.Zuazua \cite{CaZu} proved that Strichartz estimates (even with loss of derivatives) fail for metrics only $C^{0,\alpha}$ with $0\leq \alpha<1$. Our result shows a Strichartz inequality with loss of derivatives for $C^{0,1}$ metric.

Applying Theorem \ref{StrThm} for $M$ the reflection of $\Omega$ and for $G$ the reflected metric, we deduce the following theorem.

\begin{thm}
\label{MainStr}
Let $(p,q)$ be an admissible pair in dimension $d=2$ or $3$ and $I$ a finite time interval. Then there exists a constant $c(p,I)>0$ such that for any $u_0\in H_D^{\frac{3}{2p}+\epsilon}(\Omega)$ and $f \in L^1(I,H_D^{\frac{3}{2p}+\epsilon}(\Omega))$,
\begin{equation}
\label{MainStrest}
\begin{array}{rcl}
\nX{e^{it\tr}u_0}{L^p(I,L^q(\Omega))} &\leq & c(p,I) \nX{u_0}{H^{\frac{3}{2p}+\epsilon}(\Omega)},\\
\nX{\int_0^t e^{i(t-\tau)\tr} f(\tau){\rm d}\tau}{L^p(I,L^q(\Omega))} &\leq & c(p,I) \nX{f}{L^1(I, H^{\frac{3}{2p}+\epsilon}(\Omega))},
\end{array}
\end{equation}
for some $\epsilon>0$ arbitrarily small.
\end{thm}

This inequality gives us a gain of $\frac{1}{2p}-\epsilon$ derivatives with respect to the Sobolev imbedding. Compared with the Strichartz inequality obtained in the case of boundaryless Riemannian compact manifolds in \cite{BGT1} we have a supplementary loss of $\frac{1}{2p}+\epsilon$.

One could ask about the optimality of those estimates. An usual way of checking optimality is testing the estimates for solutions of the Schrödinger flow with initial data eigenfunctions of the Laplacian. This yields some $L^2\rightarrow L^q$ estimates for the eigenfunctions and we look for the optimality of those ones. We refer to some recent work of H.Smith and C.Sogge \cite{SmSo} where they prove $L^2\rightarrow L^q$ estimates for spectral clusters on regular compact domains of $\RR^d$, $d\geq 2$. Compared to those estimates, the Strichartz estimate we obtain is not optimal. Nevertheless, it has the advantage of being true for all solutions of the linear Schrödinger equation, not only those with initial data an eigenfunction. And it allows us to prove local and global existence results for the solutions of (\ref{NLS}) in dimension $d=2$.

In the proofs of local and global existence we use the $L^p(L^\infty)$ estimate of the flow in order to control the nonlinear term. We deduce it in dimension $d=2$ by combining estimates (\ref{MainStrest}) and Sobolev imbeddings.

\begin{cor}
\label{LStr}
Let $2<p<\infty$ and $d=2$.  For any $u_0\in H_D^{1}(\Omega)$ and $f \in L^1(I,H_D^{1}(\Omega))$ we have the followings inequalities
\begin{equation}
\label{TStrest}
\begin{array}{rcl}
\nX{e^{it\tr}u_0}{L^p(I,L^\infty(\Omega))} &\leq & c(p,I) \nX{u_0}{H^1(\Omega)},\\
\nX{\int_0^t e^{i(t-\tau)\tr} f(\tau){\rm d}\tau}{L^p(I,L^\infty(\Omega))} &\leq & c(p,I) \nX{f}{L^1(I, H^1(\Omega))}.
\end{array}
\end{equation} 
\end{cor}

Under this form we have a gain of $\epsilon>0$ with respect to Sobolev imbeddings (as $H_0^1(\Omega) \subset L^q$ for all $2\leq q<\infty$) by taking the average in time. However this small gain is sufficient in $\Omega\subset\RR^2$ for proving the following global existence result

\begin{thm}
\label{GEThm}
Let $\beta \in 2\mathbb{N}$, $\beta\geq 2$ and $d=2$.
For all $u_0 \in H_0^1(\Omega)$ there exists an unique solution $$u \in C(\RR, H_0^1(\Omega)) \cap L_{loc}^p(\RR,L^\infty(\Omega))$$ (for every $p>\beta$) of equation (\ref{NLS}). Moreover, for some $T>0$, the flow $u_0 \mapsto u$ is Lipschitz from $B$ bounded subset of $H_0^1(\Omega)$ to $C([-T,T],H_0^1(\Omega))$.
\end{thm}

\begin{rem} The results presented in this introduction also hold for the Schrödinger equation with Neumann boundary conditions. We shall state along the article the changes that must be done for this.
\end{rem}

\begin{rem} The Strichartz inequality also holds if $\Omega$ is the exterior of a regular bounded domain with compact boundary. We shall mention the changes that need to be made throughout the proof.
\end{rem}

This paper is organized as follows : in section \ref{NL_thm} we show how we can deduce Theorem \ref{GEThm} from Corollary \ref{LStr}. In section \ref{dbl} we present the reduction to a compact manifold endowed with a Lipschitz metric and how Theorem \ref{LStr} reads in this setting. Section \ref{Str_lips} is devoted to the proof of the Strichartz estimate.

\section{Proof of the global existence theorem}
\label{NL_thm}
Assuming the Strichartz inequality (\ref{MainStrest}), and therefore (\ref{TStrest}), we prove local existence theorem for equation (\ref{NLS}) in dimension $d=2$.
We deduce then the global existence theorem via conservation laws.

\begin{thm}
\label{LEThm}
 \textbf{(local existence theorem)} Let $\beta \in 2\mathbb{N}$, $\beta\geq 2$.
For every bounded subset $B$ of $H_0^1(\Omega)$ there exists $T>0$ such that for all $u_0 \in B$ there exists an unique solution $$u \in C([-T,T], H_0^1(\Omega)) \cap L^p([-T,T],L^\infty(\Omega))$$ (for every $p>\beta$) of equation (\ref{NLS}). Moreover, the flow $u_0 \mapsto u$ is Lipschitz from $B$ to $C([-T,T],H_0^1(\Omega))$.
\end{thm}

Note that the Lipschitz regularity of the flow was not known even in the case of cubic nonlinearity. This provides us with supplementary information about the stability of the flow under small variations of the initial data.

\begin{proof}
We denote by $X_T=C([-T,T], H_0^1(\Omega)) \cap L^p([-T,T],L^\infty(\Omega))$. This is a complete Banach space for the following norm
$$\nX{u}{X_T}={\rm max}_{|t|\leq T}\nH{u(t)}{1} + \nX{u}{L^p([-T,T],L^\infty)}.$$
We use a contraction mapping argument to deduce the existence and uniqueness of the local solution. For a fix $u_0 \in H^1$ and for $u\in X_T$ let denote by $\Phi(u)$ the functional
$$\Phi(u)(t) = e^{it\tr}u_0 -i \int_0^t e^{i(t-\tau)\tr}|u(\tau)|^\beta u(\tau) \rm{d}\tau.$$
Using the $H^1$ conservation law of the flow $e^{it\tr}$, we estimate the $H^1$ norm of $\Phi(u)$
$$
\begin{array}{rcl}
\nH{\Phi(u)(t)}{1} &\leq & \nH{u_0}{1} + \int_0^T \nH{|u(\tau)|^\beta u(\tau)}{1} \rm{d}\tau \\
&\leq & \nH{u_0}{1} +  cT^{1-\frac {\beta}{p}} \nLL{u}{p}{\infty}^\beta \nX{u}{L_T^\infty(H^1)} \leq \nH{u_0}{1} + c T^{1-\frac {\beta}{p}} \nX{u}{X_T}^{\beta +1}.
\end{array}
$$
We have used the Holder inequality to bound the $L^1$ norm of product of functions by the product of $L^p$ and $L^{p'}$ norms of functions as well as the following lemma (see e.g. \cite{AlG})

\begin{citelem} Let $s\geq 0$. Then there exists a constant $c>0$ such that for all $u,v \in H^s \cap L^\infty$ we have:
$$\nH{uv}{s} \leq c(\nH{u}{s} \nL{v}{\infty} + \nL{u}{\infty}\nH{v}{s}).$$
\end{citelem}

\noindent In order to bound the $L^p([-T,T],L^\infty(\Omega))$ norm of $\Phi(u)$ we use the Strichartz type estimate of the linear flow in this norm by the $H_0^1(\Omega)$ norm of the initial data (see estimate (\ref{LStr})).
$$
\begin{array}{rcl}
\nX{\Phi(u)}{L^p(L^\infty)} & \leq & \nLL{e^{it\tr}u_0}{p}{\infty} + \nX{\int_0^t e^{i(t-\tau)\tr}|u|^\beta u(\tau){\rm d}\tau} {L_t^p (L_x^\infty)} \\
& \leq & c \nH{u_0}{1} + c \int_0^T \nX{|u|^{\beta}u(\tau)}{H^1(\Omega)}{\rm d}\tau\\
& \leq & c \nH{u_0}{1} + cT^{1-\frac {\beta}{p}} \nX{u}{L^\infty(H^1)} \nX{u}{L^p(L^\infty)}^\beta.
\end{array}
$$
Putting together those estimates we get $\nX{\Phi(u)}{X_T} \leq c(\nH{u_0}{1} + T^{1-\frac {\beta}{p}} \nX{u}{X_T}^{\beta+1}).$ Using similar arguments we get, for $u,v\in X_T$, the following
$$
\begin{array}{rcl}
\nH{\Phi(u)(t)-\Phi(v)(t)}{1} &\leq &\int_0^T \nH{|u(\tau)|^\beta u(\tau)- |v(\tau)|^\beta v(\tau)}{1}\rm{d}\tau \\
& \leq & c \nX{u-v}{X_T}\left ( \nX{u}{X_T}^\beta + \nX{v}{X_T}^\beta \right )T^{1-\frac \beta p}
\end{array}
$$
and
$$\nX{\Phi(u) - \Phi(v)}{L_T^p L^\infty} \leq c \nX{u-v}{X_T}\left ( \nX{u}{X_T}^\beta + \nX{v}{X_T} ^\beta \right )T^{1-\frac \beta p}.$$
Let us recall that $u_0\in B$, a bounded subset of $H^1$. Then there exists $M>0$ such that for $u_0\in B$ we have $\nH{u_0}{1}\leq M$. Choose $R>0$ and $T>0$ such that $c(M+T^{1-\frac \beta p}R^{\beta+1}) < R$. For example let $R$ be $R=2cM$ and $T<cM^{-\beta \frac{p-\beta}{p}}$. This ensures that $\Phi$ maps $B(0,R;X_T)$ into $B(0,R;X_T)$. We can take $T$ even smaller and have $2c R^\beta T^{1-\frac \beta p} < 1$ and thus $\Phi$ is a contraction on $B(0,R; X_T)$. Then there exists $u\in B(0,R;X_T)$ a fixed point for $\Phi$ which therefore is the solution of equation (\ref{NLS}).
Let $u,v \in X_T$ be two solutions corresponding to two initial data $u_0,\ v_0$. We can estimate their difference uniformly in time: for all $t$ with $|t|\leq T$
$$\nX{u - v}{X_T} \leq \nH{u_0-v_0}{1} + c T^{1-\frac \beta p}\left(\nX{u}{X_T}^\beta + \nX{v}{X_T}^\beta \right) \nX{u-v}{X_T}.$$ As we have chosen $T>0$ and $R>0$ such that $2c T^{1-\frac \beta p} R^\beta <1$ we deduce the existence of a constant $C>0$ such that $\nX{u-v}{X_T} \leq C \nH{u_0-v_0}{1}$. As $\nX{u-v}{L_T^\infty H^1}\leq \nX{u-v}{X_T}$, we conclude on the Lipschitz property of the solution flow on bounded subsets of $H_0^1$.
\end{proof}

Note that this local existence theorem works for a focusing nonlinearity as well.

It is classical that when we have a Strichartz inequality, propagation of regularity holds. We give the result and a brief sketch of the proof.

\begin{prop} \textbf{(propagation of regularity)} Under same hypothesis as Theorem \ref{LEThm}, if moreover $u_0\in H^2(\Omega)$, then $u\in C([-T,T],H^2(\Omega)\cap H_0^1(\Omega))$ (same $T>0$ as in Theorem \ref{LEThm}).
\end{prop}
\begin{proof} As $u_0\in H^2(\Omega)\cap H_0^1(\Omega)\subset H_0^1(\Omega)$, we deduce from Theorem \ref{LEThm} the existence of a time $T_1>0$ such that there is a unique solution $u$ of (\ref{NLS}) in $X_{T_1}$. The same proof works for $u_0\in H^2(\Omega)\cap H_0^1(\Omega)$ and $Y_T=C([-T,T],H^2(\Omega)\cap H_0^1(\Omega))\cap L^p([-T,T],L^\infty(\Omega))$ with the norm $$\nX{u}{Y_T}={\rm max}_{|t|\leq T}\nH{u(t)}{2} + \nX{u}{L^p([-T,T],L^\infty)}.$$ Using uniqueness and $Y_T\subset X_T$, we deduce the existence of a time $0<T_2\leq T_1$ such that $u\in Y_{T_2}$. For a $T<T_2$, using the monotony of the norm $\nX{u}{X_T}$ as a function of $T$, we can establish the following inequality
$$\nX{u}{L_{T}^\infty(H^2)} \leq c_2(\nH{u_0}{2} + T^{1-\frac {\beta}{p}} \nX{u}{X_{T_1}}^{\beta}\nX{u}{L_{T}^\infty(H^2)}).$$
We take $T=T_0>0$ such that $c_2 T_0^{1-\frac{\beta}{p}} \nX{u}{X_{T_1}}<\frac 12$. This insures that $\nX{u}{L_{T_0}^\infty(H^2)}\leq 2 c_2 \nX{u_0}{H^2}$. Note that $T_0$ only depends on $\nX{u}{X_{T_1}}$. Thus we can make a bootstrap argument and conclude that $\nX{u}{L_{T_1}^\infty(H^2)}<\infty$ and thus $u\in Y_{T_1}$, for the same $T_1$ as in Theorem \ref{LEThm}.
\end{proof}

The semilinear Schrödinger equation (\ref{NLS}) has a Hamiltonian structure with gauge invariance and thus conservation laws hold for $H^2$ initial data. For $u_0\in H^1$ we deduce them by density.
\begin{prop} \textbf{(conservation laws)} The solution of (\ref{NLS}) constructed in Theorem \ref{LEThm} satisfies, for $|t|\leq T$, to 
$$\left\{ \begin{array}{l}
\int |u(t)|^2 dx = \int |u_0|^2 dx, \\
\int |\nabla u(t)|^2 + \frac{2}{\beta+2}|u(t)|^{\beta+2} dx = \int |\nabla u_0|^2 + \frac{2}{\beta+2} |u_0|^{\beta+2} dx.
\end{array}
\right.$$
\end{prop}

As a consequence, we infer the following.
\begin{thm}
\textbf{(global existence theorem)} 
The solution constructed in Theorem \ref{LEThm} extends to a global solution 
$$u \in C(\mathbb{R}, H_0^1(\Omega)) \cap L_{loc}^p(\mathbb{R},L^\infty(\Omega)).$$
\end{thm}

The proof is classical and uses the control of the $H^1$ norm by the conservation laws, as well as a bootstrap argument.

\section{Reduction to a compact manifold endowed with a Lipschitz metric}
\label{dbl}
Let $\Omega$ be a regular domain of $\RR^d$. We present here the classical mirror reflection that allows us to pass from a manifold with boundary to a boundaryless manifold. This method consists in taking a copy of the domain and glue it to the initial one by identifying the points of the boundary. In order for this to be a manifold we have to choose the coordinates carefully. Thus, taking normal coordinates at the boundary is like straightening a neighborhood of the boundary into a cylinder $\partial\Omega \times [0,1)$ and gluing the two cylinders along the boundary makes a nice smooth manifold. This can be properly done using for example tubular neighborhoods. We cite here two lemmas that can be found in \cite{Spiv}, pp. 468 and 74.

\begin{citelem}Let $\Omega$ be a regular domain of $\RR^d$, with compact boundary $\partial \Omega$. Then $\partial \Omega$ has arbitrarily small open neighborhoods in $\overline{\Omega}$ for which there are deformation retractions onto $\partial \Omega$.
\end{citelem}

The proof uses the inward pointing normal vector $\vec{n}$ and ensures the existence of a small neighborhood $U$ of $\partial \Omega$ in $\overline{\Omega}$, of a constant $\epsilon>0$ and of a diffeomorphism $\chi : U \rightarrow \partial\Omega\times [0,1)$ such that $\chi^{-1}(p,t)=p+\epsilon t \vec{n}$ for all $p\in \partial\Omega$ and $t\in [0,1)$.

Let $M=\Omega\times\{0\}\cup_{\partial\Omega}\Omega\times\{1\}$, where we identify $(p,0)$ with $(p,1)$ for $p\in \partial\Omega$. We define, for $p\in \partial\Omega$ and $t\in (-1,1)$, the mapping 
$$\tilde{\chi}^{-1}(p,t)=\left\{
\begin{array}{ccl}
(\chi^{-1}(p,t),0) & & t>0\\
p& & t=0\\
(\chi^{-1}(p,-t),1) & & t<0
\end{array} \right.$$

\begin{citelem} \textbf{\upshape{(\cite{Spiv})}}There is a unique $C^\infty$ structure on M such that $\Omega\times \{j\} \hookrightarrow M$ is $C^\infty$ and $\tilde{\chi}: U\times\{0\} \cup_{\partial\Omega} U\times\{1\} \rightarrow \partial\Omega \times (-1,1)$ is a diffeomorphism.
\end{citelem}

\textit {Note that those lemma also apply to exterior of bounded domains as long as the domains are regular and have compact boundary.}

On $M$ we define the metric $G$ induced by the new coordinates. As we have chosen coordinates in the normal direction, the metric is well defined over the boundary, its coefficients are Lipschitz in local coordinates and diagonal by blocs (no interaction between the normal and the tangent components). Moreover, $$G(r(y))=G(y),$$ where $r:M\rightarrow M,\ r(x,0)=(x,1),\ r^2=Id$ is the reflection with respect to the boundary $\partial\Omega$.

For the {\bf Dirichlet} problem we introduce the space $H_{AS}^1$ of functions of $H^1(M)$ which are anti-symmetric with respect to the boundary. Let
$$H_{AS}^1=\{v:M\rightarrow \mathbb{C},\ v\in H^1(M),\ v(y)=-v(r(y))\}.$$

\noindent Note that for $v\in H_{AS}^1$ the restriction $v_{|_{\Omega\times{0}}}$ is in $H_0^1(\Omega)$ and every function from $H_{AS}^1$ is obtained from a function of $H_0^1(\Omega)$. We shall prove the stability of $H_{AS}^1$ under the action of $e^{it\tr_G}$. 

By complex interpolation define $H_{AS}^s$ for $s\in [0,1]$ and deduce its stability under the action of $e^{it\tr_G}$. Moreover, the restriction to $\Omega$ of functions in $H_{AS}^s$ belongs to $H_D^s(\Omega)$ and vice versa. This allows us to deduce the Strichartz inequality for $e^{it\tr_D}$ on $\Omega$ from the Strichartz inequality for $e^{it\tr_G}$ on $M$. 

In section \ref{Str_lips} we give the proof of the Strichartz estimate on $(M,G)$. 

Similarly, we can define for the {\bf Neumann} problem the space $H_S^1$ of symmetric functions with respect to the boundary. This space is also stable under the action of $e^{it\tr_G}$. Thus from the Strichartz inequality on $(M,G)$ we can deduce local and global results for the Schrödinger equation (\ref{NLS}) on $\Omega$ with Neumann conditions instead of Dirichlet. Let
$$H_{S}^1=\{v:M\rightarrow \mathbb{C},\ v \in H^1(M),\ v(y)=v(r(y))\}.$$

Let us prove the stability of $H_{AS}^1$ under the action of $e^{it\tr_G}$. Let $v_0\in H_{AS}^1$ and $v(t,y)=e^{it\tr_G} v_0$. Then $v$ satisfies to $i \partial_t v(t,y)+ \tr_{G(y)} v (t,y)= 0$, $v(0)=v_0$. Let $\tilde{v}(t,y)=v(t,r(y)).$ We shall look for the equation verified by $\tilde{v}$. First note that $\tilde{v}(0)=-v_0$ and $\partial_t \tilde{v}(t,y)=\partial_t v(t,y)$. As $G$ is diagonal by blocks, having no interactions between the normal and tangent components, so is $G^{-1}$. Thus in $\tr_{G(y)}$ there is no crossed term. Consequently $\tr_{G(r(y))}\tilde{v}(t,y) = \tr_{G(y)}v(t,y)$. We see thus that $\tilde{v}$ satisfies to the linear Schrödinger equation with initial data $-v_0(y)$. But $-v(t,y)$ satisfies the same equations. By uniqueness we conclude that $$v(t,r(y))=-v(t,y).$$

We are now able to prove the following proposition.
\begin{prop}
Theorem \ref{StrThm} implies Theorem \ref{MainStr}.
\end{prop}

\begin{proof} Let $M$ be the reflection of $\Omega$ and $G$ the reflected metric. Consider $u_0\in H_D^{\frac{3}{2p}+\epsilon}(\Omega)$. Let $v_0:M\rightarrow \mathbb{C}$ be defined as follows : for $y\in \Omega$, let $v_0((y,0))=u_0(y)$ and $v_0((y,1))=-u_0(y)$. As seen previously $v_0\in H_{AS}^{\frac{3}{2p}+\epsilon}\subset H^{\frac{3}{2p}+\epsilon}(M)$. Moreover $\nX{v_0}{H^{\frac{3}{2p}+\epsilon}(M)}^2=2\nX{u_0}{H^{\frac{3}{2p}+\epsilon}(\Omega)}^2$. From the stability of the $H_{AS}^{\frac{3}{2p}+\epsilon}$ under the action of $e^{it\tr_G}$ and the uniqueness of the linear flow we conclude that $e^{it\tr_G}v_0{|_{ \Omega\times\{0\}}} = e^{it\tr}u_0$. This leads us to 
$$\nX{e^{it\tr}u_0}{L^p(I,L^q(\Omega))} \leq c(p,I) \nX{u_0}{H^{\frac{3}{2p}+\epsilon}(\Omega)},$$
which is the first estimate (\ref{MainStrest}) in Theorem \ref{MainStr}. Estimate in the nonhomogeneous form is obtained classically by means of Minkowski inequality from the homogeneous estimate (see e.g. \cite{BGT1}). 
\end{proof}

In the next section we prove Theorem \ref{StrThm}.

\section{Strichartz inequality for the Schrödinger operator associated to a Lipschitz metric}
\label{Str_lips}
Let $M$ be a $C^\infty$ compact manifold (or flat outside a compact set) endowed with a metric
whose coefficients are Lipschitz. We want to study the behavior of the Schrödinger flow in the $L_t^p(L_x^q)$ norm and for doing so we translate the equation in local coordinates of $\RR^d$. Having a Schrödinger equation we pass in semiclassical time coordinates and study frequency localized initial data restricted to a coordinate chart (in this way the solution remains essentially localized in the open chart on a very short time that depends on the frequency, as we shall see). We use a partition of unity to recover the behavior of the solution on the whole manifold.

\subsection{Preliminaries}
\label{partition}

In the case $M$ compact manifold, let $(U_j,\kappa_j)_{j\in J}$ be a finite covering with open charts. Let $(\chi_j)_{j\in J} : M \rightarrow [0,1]$ be a partition of unity subordinated to the covering $(U_j)_{j\in J}$. For all $j \in J$ let $\tilde{\chi}_j : M \rightarrow [0,1]$ be a $C^\infty$ function such that $\tilde\chi_j=1$ on the support of $\chi_j$ and the support of $\tilde\chi_j$ is contained in $U_j$ .

The coordinate map $\kappa_j : U_j \subset M \rightarrow V_j \subset \RR^d$ transports the functions $\chi_j$ and $\tilde\chi_j$ onto the functions $\chi_j^1(y)=\chi_j(\kappa_j^{-1}(y))$ and $\chi_j^2(y)=\tilde\chi_j(\kappa_j^{-1}(y))$.

In the case $M$ flat outside a compact set, let $(U_j,\kappa_j)_{j\in J}$ be a covering of the area of $M$ where $G\neq \mathbb{I}d$. This area is compact, so we can choose $J$ of finite cardinal. We have $M=\cup_{j\in J}U_j \cup U_{1,\infty}\cup U_{2,\infty}$, where $U_{1,\infty}$ and $U_{2,\infty}$ are two disjoint neighborhood of $\infty$, diffeomorphe to $\RR^d \backslash \bar{B}$. Let $(\chi_j)_{j\in J},\ \chi_{1,\infty},\ \chi_{2,\infty} : M \rightarrow [0,1]$ be a partition of unity subordinated to the previous covering. We estimate $e^{it\tr_G}u_0$ on $\cup_{j\in J}U_j$ exactly as we do for the compact manifold. Knowing that $G = \mathbb{I}d$ on $U_{\infty}$ simplifies the analysis of the spectrally truncated flow near infinity. 

We prepare the frequency decomposition. Let $\varphi_0 \in C^\infty(\RR^d)$ be supported in a ball centered at origin and $\varphi \in C^\infty(\RR^d)$ be supported in an annulus such that for all $\lambda \in \RR^d$
\begin{equation}
\label{LPdecomp}
\varphi_0(\lambda) + \sum_{k \in \mathbb{N}} \varphi(2^{-k} \lambda)=1.
\end{equation}

We define a family of spectral truncations : for $f\in C^\infty(M)$ and $h\in (0,1)$ let
\begin{equation}
J_h f=\sum_{j\in J} (\kappa_j)^* \left( \chi_j^2 \varphi(hD)(\kappa_j^{-1})^* (\chi_j f) \right)
\end{equation} and
\begin{equation}
J_0 f=\sum_{j\in J} (\kappa_j)^* \left( \chi_j^2 \varphi_0(D)(\kappa_j^{-1})^* (\chi_j f) \right),
\end{equation}
where $^*$ denotes the usual pullback operation. We can rewrite $J_h$ as follows
$$J_h f(x)= \sum_{j\in J} \tilde{\chi_j}(x) \varphi(hD) \left(\chi_j(\kappa_j^{-1}) f(\kappa_j^{-1}) \right)(\kappa_j(x)).$$
If we denote by $\rho$ and $\rho_0$ the inverse Fourier transform of $\varphi$ and $\varphi_0$ respectively and if we set $f_j=\chi_j(\kappa_j^{-1})f(\kappa_j^{-1})$, then
$$\varphi(hD)f_j(\kappa_j(x))= \frac{1}{h^d}\int_{\RR^d}\rho \left(\frac {\kappa_j(x)-z}{h}\right)f_j(z){\rm d}z.$$
From relation (\ref{LPdecomp}) we deduce that for all $x\in \RR^d$ and for $v$ a function on $\RR^d$ :
$$\int_{\RR^d} \rho_0(x-y) v(y){\rm d}y + \sum_{k=0}^\infty 2^{kd}\int_{\RR^d}\rho \left( 2^k(x-y) \right) v(y) {\rm d}y = v(x).$$
We obtain thus
\begin{equation}
\label{2LPdecomp}
J_0 f(x) + \sum_{k=0}^\infty J_{2^{-k}} f(x) = \sum_{j\in J} \tilde{\chi_j}(x) \chi_j(\kappa_j^{-1}(\kappa_j(x))) f(\kappa_j^{-1}(\kappa_j(x)))=f(x).
\end{equation}

Note that in the case $M$ flat outside a compact set, we have to modify $J_h$ such that it takes into account the influence of the spectral truncation on the chart near $\infty$. Let
\begin{equation}
\label{LPdecompFlat}
F_{\infty} f =\tilde{\chi}_{\infty} \varphi(hD) \chi_{\infty}f(x).
\end{equation}
Then for $J_{h,\infty}=J_h + F_{1,\infty}+F_{2,\infty}$ we have an identity similar to (\ref{2LPdecomp}).

We study the semiclassical Schrödinger equation with initial data $J_h u_0$ and then we recover the behavior of the linear flow thanks to identity (\ref{2LPdecomp}). We introduce the semiclassical time $s$ by $w(s,x)=v(hs,x)$. If $v$ is a solution of the equation
$$\left\{
\begin{array}{rcl}
i \partial_t v + \tr_G v & = & 0\\
v_{|_{t=0}} & = & J_h u_0
\end{array}\right.
$$
on a time interval $I$, then $w$ is solution of the following semiclassical equation on $h^{-1}I$
\begin{equation}
\label{Sscl} \left\{
\begin{array}{rcl}
i h\partial_s w + h^2 \tr_G w & = & 0\\
w_{|{s=0}} & = & J_h u_0.
\end{array}\right.
\end{equation}

The classical way of proving Strichartz inequalities is to use the $TT^*$ method (here $^*$ stays for adjoint) starting from
a $L^2$ conservation norm and a $L^1 \rightarrow L^\infty$
dispersive estimate (\cite{GiVe}). In the case under study, the dispersive estimate can be
obtained by combining the WKB approximation (as in \cite{BGT1}) and a stationary phase
type lemma. In order to use this strategy we need more regularity on
the coefficients of the metric. Using an idea from \cite{BaCh} (see also \cite{Ta}), we regularize them at some frequency $h^{-\alpha}$, where $0<\alpha<1$ is  a parameter that will be fixed in the end. We treat the remainder term as a source term like in \cite{BaCh}.

Let $\psi$ be a $C_0^\infty(\RR^d)$ radially symmetric function with $\psi\equiv 1$ near 0. We define the regularized metric $G_h$ as follows
\begin{equation}
\label{regmetricdef}
G_h =\sum_{j\in J} (\kappa_j)^* \left( \chi_j^2 \psi(h^\alpha D)(\kappa_j^{-1})^* (\chi_j G) \right).
\end{equation}
The transformation of $G$ into $G_h$ does not spoil the symmetry. Note also that $G_h$ converges uniformly in $x$ to $G$, and thus, for $h$ sufficiently small, $G_h$ is positive definite. Therefore, $G_h$ is still a metric. Then equation (\ref{Sscl}) is equivalent to
$$\left\{
\begin{array}{rcl}
i h\partial_s w + h^2 \tr_{G_h} w & = & h^2 (\tr_{G_h} - \tr_G)w\\
w_{|{s=0}} & = & J_h u_0.
\end{array}\right.
$$
When writing $J_h$ in local coordinates we see it as a finite sum of expressions as
\begin{equation}
\label{defFh}
F_h f(x) = \tilde{\chi}(x) \varphi(hD)(\chi f)(x)= \frac{1}{h^d} \int_{\RR^d} \tilde{\chi}(x)\rho \left( \frac{x-y}{h} \right)f(y){\rm d}y,
\end{equation}
where $x\in \RR^d$, $f:\RR^d \rightarrow \RR$, $\chi$ and $\tilde{\chi}$ are compactly supported , $0\leq \chi,\ \tilde{\chi} \leq 1$ and $\tilde{\chi}\equiv 1$ on the support of $\chi$. The function $\varphi$ is $C^\infty$ supported in an annulus. We study the following equation in local coordinates
\begin{equation}
\label{Shcoord}
\left\{
\begin{array}{rcl}
i h\partial_s w + h^2 \tr_{G_h} w & = & 0\\
w_{|{s=0}} & = & F_h u_0.
\end{array}\right.
\end{equation}

The plan of the proof is the following :
\begin{itemize}
\label{plan}
	\item{construct an approximate solution for (\ref{Shcoord}) by the WKB method and prove the dispersion estimate on a small interval of time $I_h$.  This solution remains supported in the chart domain so we can extend it as a function onto the manifold.}
	\item{obtain a Strichartz inequality for the spectrally truncated flow $J_h^* e^{ihs\tr_{G_h}}$ on $I_h$, where $J_h^*$ denotes the $L^2$ adjoint of $J_h$.}
	\item{estimate the difference between the regularized flow and the initial flow in the $L^p(L^q)$ norm on $I_h$.}
	\item{obtain the Strichartz inequality for $e^{it\tr_G}$ on a fixed time interval.}
\end{itemize}

The analysis of $F_\infty^* e^{it\tr_G}u_0$ in $L^p(L^q)$ norm on a small interval of time can be done using the classical Strichartz estimate (see proof of Proposition \ref{PFIhStr}).

\subsection{Estimates on the regularized metric and preliminary commutator lemmas}
The metric $G:M\rightarrow M_d(\RR)$ is symmetric, positive definite and Lipschitz : there exist $c,C,c_1 >0$ such that for all $x\in M$
$$c{\mathbb{I} d} \leq G(x) \leq C {\mathbb{I} d},\ |\partial G|\leq c_1,$$
where we have denoted by $\partial G$ the derivatives of the metric in a system of coordinates. Using the expression (\ref{regmetricdef}), one can easily prove the following estimates

\begin{prop}
\label{regmetric} The regularized metric $G_h$ is a $C^\infty$ function
that verifies, in a system of coordinates, the followings : there exists $c, C>0$ and $c_\gamma >0$ for all
$\gamma \in \mathbb{N}^d$ such that for all $x\in M$
$$c {\mathbb I}d \leq G_h(x) \leq C {\mathbb I}d,\ |\partial^\gamma G_h(x)| \leq c_\gamma h^{-\alpha \max (|\gamma|-1 , 0)}.$$
\end{prop}

Next, we present a collection of useful lemmas about the action of operators $F_h$ defined in (\ref{defFh}).

\begin{lem} There exists a constant $C>0$ such that, for all $1\leq p \leq \infty$, $F_h$ is bounded from $L^p$ to $L^p$
$$\nX{F_h}{L^p \rightarrow L^p} \leq C.$$
\end{lem}
\begin{proof}
If we denote by $f_1=\chi f$ then the boundedness of $\chi$ ensures $\nL{f_1}{p}\leq \nL{f}{p}$. Thus, the result follows from the classical estimate $\nX{\varphi(hD)}{L^p\rightarrow L^p}\leq C$.
\end{proof}

\begin{lem}
\label{commJhGh}
There exist constants $c_1>0$ and $c_2>0$ such that the commutator $[F_h,\tr_{G_h}]=F_h \tr_{G_h} - \tr_{G_h} F_h$ is bounded from $L^2$ to $L^2$ of norm $\frac{c_1}{h}$ and from $H^1$ to $L^2$ of norm $c_2$:
$$\nX{[F_h,\tr_{G_h}]}{L^2\rightarrow L^2} \leq \frac{c_1}{h}$$
and
$$\nX{[F_h,\tr_{G_h}]}{H^1\rightarrow L^2} \leq c_2.$$
\end{lem}

\noindent We shall use the following

\begin{Schur} For $T$ a kernel operator, $Tf(x)=\int_{\RR^d}K(x,y)f(y) {\rm d}y$, if $$\max \left( \int_{\RR^d} |K(x,y)| {\rm d}y,\ \int_{\RR^d}|K(x,y)| {\rm d}x \right) \leq c,$$ then for all $1 \leq p\leq \infty$ we have $T:L^p(\RR^d)\rightarrow L^p(\RR^d)$ and $\nX{T}{L^p\rightarrow L^p}\leq c$.
\end{Schur}

\begin{proof} \emph{of Lemma \ref{commJhGh}.}
We first prove the $L^2\rightarrow L^2$ estimate. 
We write the commutator $[F_h,\tr_{G_h}]$ as a convolution kernel operator $[F_h,\tr_{G_h}]=\frac{1}{h^d}\int_{\RR^d}k_{1}(x,\frac{x-y}{h})f(y){\rm d}y$ by doing integration by parts. We arrange the terms in $k_1$ according to the order of derivatives on $\rho$. We estimate the coefficients in $L^\infty$ norm. The coefficients of $\rho$ must have 2 derivatives on $\chi$ or $G_h$. The biggest among them is the one where both derivatives bear on $\frac{1}{\sqrt{\det G_h(y)}}$. By Proposition \ref{regmetric}, this term is of order $h^{-\alpha}(<h^{-1})$. All other coefficients of $\rho$ are bounded. The coefficients of $\frac 1h \partial_j \rho(\frac{x-y}{h})$ have one derivative on $G_h$ or $\chi$ and thus are bounded. The coefficient of $\frac{1}{h^2}\partial_i\partial_j\rho(\frac{x-y}{h})$ is $$\tilde{\chi}(x)(G_h^{i,j}(y)-G_h^{i,j}(x))\chi(y).$$ It is of order $|x-y|$ and if we denote by $\rho_1^{i,j}(\frac{x-y}{h})=\frac{x-y}{h}\partial_i \partial_j \rho(\frac{x-y}{h})$ then $\frac{1}{h^d}\rho_1(\frac{x-y}{h})$ satisfies the conditions from Schur's lemma. We conclude that the $L^2\rightarrow L^2$ norm of the commutator is of order $h^{-1}$.

For the $H^1\rightarrow L^2$ estimate we write the commutator as a convolution kernel operator as follows
$$[F_h,\tr_{G_h}] f(x)= \frac{1}{h^d}\int_{\RR^d}k_{0}\left( x,\frac{x-y}{h}\right) f(y){\rm d}y+ \frac{1}{h^d}\int_{\RR^d}k_{2}\left(x,\frac{x-y}{h}\right) \nabla f(y){\rm d}y.$$ 
Indeed, using the obvious identity $\partial_{x_j}\left( \rho\left(\frac{x-y}{h}\right)\right) = - \partial_{y_j}\left( \rho \left( \frac{x-y}{h} \right) \right),$
we can make an integration by parts and obtain both terms in $f(y)$ and in $\nabla f(y)$. We are doing this as follows : if no derivative bears on $f$ but there is one on $\rho$, we proceed to the integration by parts.

Thus, $k_0(x,\frac{x-y}{h})$  contains no derivative of $\rho$ and therefore the operator associated to $k_0$ is bounded from $L^2$ to $L^2$. As above, we arrange the terms in $k_2$ following the order of derivatives on $\rho$. As we have at most one derivative that acts on each term, the coefficient of $\rho(\frac{x-y}{h})$ is bounded. As for the coefficient of $\frac{1}{h} \partial_i \rho(\frac{x-y}{h})$, it equals $\tilde{\chi}(x)\left( G_h^{i,j}(y) - G_h^{i,j}(x)\right) \chi(y)$ and as above we deduce the boundedness of the commutator from $H^1$ to $L^2$.   
\end{proof}

As one may not apply two derivatives on $G(x)$, the similar statement for $[F_h,\tr_G]$ only holds for the $H^1 \rightarrow L^2$ norm, namely :
$$\nX{[F_h,\tr_G]}{H^1\rightarrow L^2} \leq c.$$

\begin{lem}
\label{diffGGh}
There exists a constant $c>0$ such that the operator $F_h(\tr_{G_h}-\tr_G)$ is bounded from $H^1$ to $L^2$ with norm $ch^{\alpha -1}$,
$$\nX{F_h(\tr_{G_h}-\tr_G)}{H^1 \rightarrow L^2} \leq c h^{\alpha -1}.$$
\end{lem}
\begin{proof}
We write $F_h(\tr_G -\tr_{G_h})f$ as a convolution kernel operator that acts on $\nabla f$. We do a similar analysis of the kernel of $F_h(\tr_{G_h}-\tr_G)$ with the one done in the proof of Lemma \ref{commJhGh}. The coefficient of $\rho(\frac{x-y}{h})$ is bounded since it contains one derivative of $G$, $G_h$ or $\chi$. The coefficient of $\frac 1h \partial_i \rho(\frac{x-y}{h})$ is $\tilde{\chi}(x)\left( G_h^{i,j}(y)-G^{i,j}(y)\right)\chi(y)$. Let us recall that $G_h=\psi(h^\alpha D)G$ and $\psi(0)=1$. Thus, $$\nX{G_h-G}{L^\infty}\leq ch^\alpha.$$ The result follows from Schur's lemma. 
\end{proof}

Let $\tilde{\varphi}$ be a $C^\infty$ function supported in an annulus such that $\tilde{\varphi}=1$ on a neighborhood of the support of $\varphi$. We define $\tilde{F}_h$ just like $F_h$, replacing $\varphi$ par $\tilde{\varphi}$ (see (\ref{defFh})). We denote by $$T_h =\tilde{F}_h F_h - F_h.$$ The following lemma states that the action of $\tilde{F}_h$ on $F_h$ and $[F_h, \tr_{G_h}]$ is close to identity in $L^p\rightarrow L^p$ and $L^2 \rightarrow L^2$ norm respectively.

\begin{lem} 
\label{Rlemma}
For all $N\in \mathbb{N}$ and $p\geq 2$, the following inequalities hold
\begin{equation}
\label{Rtilde}
\nX{[F_h,\tr_{G_h}]- [F_h,\tr_{G_h}]\tilde{F}_h}{L^2\rightarrow L^2} \leq c_N h^N
\end{equation}
and
\begin{equation}
\label{RRest}
\nX{T_h}{L^p\rightarrow L^p}\leq c_N h^N.
\end{equation}
\end{lem}

\begin{proof}
As in the proof of Lemma \ref{commJhGh} we write 
$[F_h,\tr_{G_h}]\tilde{F}_h f(x) = \frac{1}{h^d}\int_{\RR^d} \tilde{k}_{1} (x,\frac{x-y}{h}) f(y) {\rm d}y$, where
$\tilde{k}_1\left( x,\frac{x-y}{h}\right) =\frac{1}{h^d}\int_{\RR^d} k_1 \left(x,\frac{x-r}{h}\right) \tilde{\chi}(r) \tilde{\rho} \left( \frac{r-y}{h} \right) \chi(y){\rm d}r.$
Using $(D^\gamma \varphi) \tilde{\varphi}= D^\gamma \varphi$ and basic properties of convolution and Fourier transform, we obtain the following identity 
\begin{equation}
\label{Fid}
\int_{\RR^d} z^\gamma \rho(z) \tilde{\rho}\left( \frac{x-y}{h}-z \right){\rm d}z= \left( \frac{x-y}{h} \right) ^\gamma \rho \left( \frac{x-y}{h} \right).
\end{equation} 

We shall use identity (\ref{Fid}) to show that
\begin{equation}
\label{SchrId}
\tilde{k}_1 \left( x,\frac{x-y}{h} \right)=k_1 \left( x,\frac{x-y}{h}\right) +h^N R_N (x,y),
\end{equation}
for all $N\in \mathbb{N}$ and such that $R_N$ satisfies conditions of Schur's lemma with a constant independent of $h$.

The kernel $k_1(x,\frac{x-r}{h})$ is a sum of terms as $\rho(\frac{x-r}{h})c_0(x,r)$ and $\frac 1h \rho_1(\frac{x-r}{h}) c_1(x,r)$, where $c_0$ and $c_1$ are factors of $G_h$ and $\chi$ (as well as their derivatives up to order 2) considered in $x$ or $r$. Here $\rho_1(\frac{x-r}{h})$ can be either $\partial_j \rho(\frac{x-r}{h})$ or $(\frac{x-r}{h})\partial_i\partial_j \rho(\frac{x-r}{h})$. We perform a Taylor expansion in $r$ of factors from $c_0$ and $c_1$ and express them in $x$. Thus
$$c_0(x,r)=\sum_{n= 0}^{N_0} \sum_{|\gamma|=n}T_{0,\gamma}(x) \left(\frac{x-r}{h}\right)^\gamma h^{|\gamma|} + R_{0,N_0}(x,r).$$
Note that $T_{0,\gamma}$ may contain derivatives of $G_h$ up to order $|\gamma|+2$. Therefore $\nL{T_{0,\gamma}}{\infty} \leq c h^{-\alpha(\gamma+1)}$. The remainder term $R_{0,N_0}(x,r)$ is of order $O(h^{-\alpha (N_0+2)} |x-r|^{N_0+1}).$ 
We will use this Taylor expansion in both directions.
First we use it to expand $c_0$ as a sum. 
By the change of variable $r=x-hz$, identity (\ref{Fid}) and the Taylor expansion from the right hand side to the left hand side, we obtain $$\frac{1}{h^d} \int \rho\left(\frac{x-r}{h} \right) c_0(x,r)\tilde{\rho} \left(\frac{r-y}{h} \right) {\rm d}r = c_0(x,y) \rho\left(\frac{x-y}{h} \right) + I_{1,N_0}(x,y) + I_{2,N_0}(x,y).$$
Here $I_{1,N_0}(x,y)$ denotes the integral with the remainder term from the Taylor expansion $I_{1,N_0}(x,y)=\frac{1}{h^d}\int \rho\left(\frac{x-r}{h}\right) R_{0,N_0} (x,r) \tilde{\rho}\left(\frac{r-y}{h}\right) {\rm d}r$ and $I_{2,N_0}(x,y)=-R_{0,N_0}(x,y) \rho\left(\frac{x-y}{h} \right)$.

We analyze the $I_{1,N_0}(x,y)$ term. For all $x,y\in \RR^d$,
$$|I_{1,N_0}(x,y)|\leq c h^{(N_0+1) (1-\alpha)-\alpha+d} \frac{1}{h^d} \int |\rho(z)| |z|^{N_0+1} |\tilde{\rho}|\left(\frac{x-y}{h}-z \right) {\rm d}z.$$ 
Thus, Schur's lemma applies for the kernel $I_{1,N_0}(x,y)$ with a constant $c h^{(N_0+1) (1-\alpha)-\alpha+d}$.
Similarly, Schur's lemma applies for $I_{2,N_0}(x,y)$ with a constant $ch^{(N_0+1) (1-\alpha)-\alpha+d}$.
We treat the $\rho_1$ term in a similar manner.

As $0<\alpha<1$, for all $N\in\mathbb{N}$ there exist a $N_0$ such that $(N_0+1) (1-\alpha)-\alpha >N$. If we denote by $R_N(x,y)=h^{-N} \left( I_{1,N_0}(x,y) +I_{2,N_0}(x,y)+... \right)$ (the $...$ stand for the remainder terms in $\rho_1$), then $R_N$ satisfies (\ref{SchrId}). As $h^N R_N(x,y)$ is the kernel of $[F_h,\tr_{G_h}]- [F_h,\tr_{G_h}]\tilde{F}_h$, inequality (\ref{Rtilde}) follows from Schur's lemma.

We now pass to the proof of (\ref{RRest}). The method of proof is very similar. Using that $\chi \tilde{\chi}=\chi$ we can write $\tilde{F_h}F_h$ as a kernel operator $\tilde{F_h}F_h f(x) = \frac{1}{h^d} \tilde{\chi}(x) \int k(x,y)\chi(y)f(y){\rm d}y,$ where $$k(x,y)=\frac{1}{h^d}\int \tilde{\rho} \left( \frac{x-r}{h} \right) \rho \left( \frac{r-y}{h} \right) \chi(r) {\rm d}r.$$ As above, using the change of variable $z=y+hz$, making a Taylor expansion of $\chi$ in $y$ and using identity (\ref{Fid}) we conclude that the kernel of $\tilde{F}_h F_h$ equals the kernel of $F_h$ plus some remainder terms. The result follows from the analysis of the remainder terms and Schur's lemma as above.
\end{proof}

\subsection{Construction and estimate of the ansatz}
\label{ansatz}

We shall construct and estimate an approximate solution on a bounded open chart. The proof need to be slightly modified to apply also for a neighborhood of $\infty$, but we shall not use it here. Let us recall the notations for the truncation in space coordinates as introduced in section \ref{partition}. We have $U\in \RR^d$ a bounded open chart. Let $\chi$ and $\tilde{\chi}$ be $C^\infty$ functions supported in $U$ such that $\tilde{\chi}\equiv 1$ on a neighborhood of the support of $\chi$.

The WKB method consists in searching for an approximate solution of equation
(\ref{Shcoord}) that decomposes as :
\begin{equation}
\label{WKBap} w_N^{ap}(s,x) = \int_{\RR^d} e^{i
\frac{\Phi(s,x,\xi)}{h}}\sum_{j=0}^N h^j a_j(s,x,\xi)
\widehat{v}_0 \left( \frac{\xi}{h}\right) \frac{d\xi}{(2\pi h)^d},
\end{equation}
with $\Phi(0,x,\xi)=x\cdot \xi$, $a_0(0,x,\xi)=\tilde{\chi}(x)\varphi(\xi)$ and
for $j\geq 1$, $a_j(0,x,\xi)=0$. We have denoted by $v_0= \chi u_0$. Thus, by the inverse Fourier
transform, $w_N^{ap}(0,x) = F_h u_0(x)$. We want $w_N^{ap}$ to be
close to the solution of (\ref{Shcoord}). In other words we want to
find $r_{h,N}$ small (in a sense that will be stated further)
such that:
\begin{equation}
\label{appeq} \left\{ \begin{array}{rcl}
i h\partial_s w_N^{ap} + h^2 \tr_{G_h} w_N^{ap} & = & r_{h,N}\\
{w_N^{ap}}_{|{s=0}} & = & F_h u_0.
\end{array}\right.
\end{equation}

If we introduce formally (\ref{WKBap}) into the equation
(\ref{Shcoord}) we see that $\Phi$ should satisfy the following
Hamilton-Jacobi equation:
\begin{equation}
\label{HJ} \left\{
\begin{array}{rcl}
\partial_s \Phi + G_h^{l,m} \partial_{x_l} \Phi \partial_{x_m} \Phi & = & 0\\
\Phi_{|{s=0}} & = & x\cdot \xi
\end{array}\right.
\end{equation}
and $a_0$ should satisfy the linear transport equation:
\begin{equation}
\label{tr0} \left\{
\begin{array}{rcl}
\partial_s a_0 + \tr_{G_h} \Phi\cdot a_0 + 2 G_h^{l,m} \partial_{x_l} \Phi \partial_{x_m} a_0 & = & 0\\
{a_0}_{|{s=0}} & = & \tilde{\chi}(x) \varphi(\xi)
\end{array}\right.
\end{equation}
while for $j \geq 1$, the $a_j$ should satisfy the nonhomogeneous
transport equation (we consider $i \tr_{G_h} a_{j-1}$ as a source term)
\begin{equation}
\label{trj} \left\{
\begin{array}{rcl}
\partial_s a_j + \tr_{G_h} \Phi\cdot a_j + 2 G_h^{l,m} \partial_{x_l} \Phi \partial_{x_m} a_j & = & i \tr_{G_h} a_{j-1}\\
{a_k}_{|{s=0}} & = & 0
\end{array}\right.
\end{equation}
Note that the functions $\Phi$ and $a_j$ depend on $h$ and this dependence will be quantified in Proposition \ref{phf_estimates}.
We recall a transport lemma that will be used in the following proofs.

\begin{lem}
\label{trlemma}
Let $f_l : \RR \times \RR^d \rightarrow \RR$ a sequence of bounded $C^1$ functions, for $1\leq l \leq d$, and $b : \RR \times \RR^d \rightarrow \RR$ a $C^1$ function such that there exists $M>0$ that bounds $|b(s)| \leq M$ for all $s\in \RR$. For $u_0 : \RR^d \rightarrow \RR$ the solution $u$ of the transport equation
$$\partial_s u +\sum f_l \partial_{x_l}u + bu=0,\ u_{|_{s=0}}=u_0,$$
satisfies
$$\nLx{u(s)}{\infty} \leq e^{M|s|} \nL{u_0}{\infty}.$$
Under the same assumptions on $f_l$ and $b$ and $F:\RR\times \RR^d \rightarrow \RR$, the solution $v$ of the nonhomogeneous transport equation
$$\partial_s v +\sum f_l \partial_{x_l}v + bv=F,\ v_{|_{s=0}}=v_0,$$
satisfies the following estimate:
$$\nLx{v(s)}{\infty} \leq e^{M|s|} \nL{v_0}{\infty} + s e^{M|s|}\nX{F}{L_{x,s}^\infty}.$$
\end{lem}
The proof is classical, using the methods of characteristics to transform the transport equation into a system of ODEs. The second part uses similar arguments combined with Gronwall lemma.

\begin{prop}
\label{time_exist} Let $R>0$ such that ${\rm supp}\ \!\varphi \subset
B(0,R)$. Then there exists $c>0$ and $\Phi \in
C^{\infty}([-ch^\alpha, c h^\alpha] \times \RR^d \times B(0,R))$
solution of the Hamilton-Jacobi equation (\ref{HJ}). There exist
$(a_j)_{j\in \mathbb{N}}$ a sequence of functions in
$C^{\infty}([-ch^\alpha, c h^\alpha] \times \RR^d \times B(0,R))$
solutions to the transport equations (\ref{tr0}) and (\ref{trj}). Moreover, the support of $a_j(s,\cdot,\xi)$ is included in $U$ (and therefore compact) for all $|s|\leq ch^\alpha$ and $\xi \in B(0,R)$.
\end{prop}

\begin{proof}
We solve the Hamilton-Jacobi equation by the method of
characteristics. For a fixed $\xi \in B(0,R)$, the symbol of the Hamiltonian is
$p(x,\eta)=-G_h^{l,m}(x) \eta_l \eta_m$. 
If we denote $\psi(s)=(y(s),\eta(s))$ with $y(s), \eta(s):
\RR^d \rightarrow \RR^d$ then the couple $(y(s),\eta(s))$ verifies the Hamiltonian system given by $p(x,\eta)$. Moreover, we impose $x\in \RR^d \mapsto y(s,x)\in \RR^d$ to be a diffeomorphism for all $\xi \in B(0,R)$ and $s\in [-S,S]$, .
We conclude by the Cauchy Lipschitz theorem the local
existence and uniqueness of smooth solutions.
As the Hamiltonian is constant on the characteristics, for all $s\in [-S,S]$,
$$-G_h^{l,m}(y(s,x,\xi)) \eta_l(s,x,\xi) \eta_m(s,x,\xi)= -G_h^{l,m}(x)\xi_l \xi_m.$$ From the equivalence of the metric $G_h$ and the metric $G$ (see Proposition \ref{regmetric}), we deduce the existence of two
constants $c,C>0$ such that for all $s,x,\xi$
\begin{equation}
\label{Phone}
c \leq |\eta|(s,x,\xi) \leq C.
\end{equation}
We have to find a time length $S>0$ such that, for all $s\in [-S,S],$ $x\in\RR^d \mapsto y(s,x)$ is a diffeomorphism of $\RR^d$. We consider the equation verified by
\begin{equation*}
J(s)={\rm det} \left( (\partial_{x_k} y_l(s))_{l,k} \right).
\end{equation*}

As $y(0,x)=x$ we have $J(0)=1$. In order to find the equation verified by $J$, we differentiate the characteristic system following $x_k$. We obtain
$$\partial_s J(s)= \sum_k {\rm det}\left(\partial_{x_1}y, \ldots,
\partial_{x_{k-1}}y, B_{11}\partial_{x_k}y + B_{12}\partial_{x_k}\eta, \partial _{x_{k+1}}y, \ldots, \partial_{x_d} y \right),$$ where $B_{11}=( -2 \partial_{x_r} G_h^{j,m}(y)\eta_m )_{j,r}$ and $B_{12}=( -2 G_h^{j,r} (y))_{j,r}$ are $d\times d$ matrices. We obtain
$$\dot{J}(s)={\rm tr}(B_{11}) J(s) + f(s),$$ where ${\rm tr}(B_{11})$ denotes the trace of $B_{11}$ and $f$ gathers all the terms that contain $\eta$. From estimate (\ref{Phone}) combined with estimates on the regularized metric we deduce $|{\rm tr}(B_{11})|\leq c$. Using Duhamel formula we get 
$$|J(s)| \geq e^{\int_0^{|s|} {\rm tr} B_{11}(r) dr} - \int_0^{|s|} e^{\int_r^{|s|} {\rm tr}
B_{11}(\tau)d\tau}|f(r)|dr.$$ 
We are looking for a $S>0$ such that, for $|s|<S$, the right hand side is strictly positive. We shall start by estimating, for all $s\in \RR$ and $x\in \RR^d$, the force term $f$. Applying Gronwall lemma to the linear system obtained differentiating the Hamiltonian system following $x_k$, we obtain $|f(s)|\leq c e^{d|s|(c+ch^{-\alpha})}$. It suffices to have $|s| \leq ch^{\alpha}$ in order to have $|f(s)|$ bounded for all $h>0$. Thus, by
taking $S=c h^\alpha$ eventually with a smaller constant $c>0$, we get
$J(s) \geq \delta >0$ for all $-S< s <S$ and therefore $x\mapsto y(s,x)$ is a diffeomorphism of $\RR^d$. 

By the method of characteristics we know $\nabla_x \Phi(s,y(s,x,\xi),\xi)=\eta(s,x,\xi)$. Inverting $x\mapsto y(s,x)$ for $|s|\leq ch^\alpha$ we obtain the announced properties for $\Phi$. Moreover, from (\ref{Phone}) we deduce
\begin{equation}
\label{important}
c \leq \nL{\nabla_x\Phi(s,x,\xi)}{\infty} \leq C.
\end{equation}
Using $a_{0{|_{s=0}}}=\tilde{\chi}(x) \varphi(\xi)$, the boundedness and the uniformity of the speed of propagation, we can take the time length $S=ch^\alpha$, with $c>0$ being chosen eventually smaller, such that $x\mapsto a_0(s,x,\xi)$ is supported in $U$ for all $|s|\leq S$ and $\xi\in B(0,R)$. 

Moreover, the equations verified by $a_j$, for $j\geq 1$, are nonhomogeneous linear equations (\ref{trj}) with initial data $0$ and source term $i\tr_{G_h}a_{j-1}$. Consequently, the support of $a_j$ is the same as the support of $a_{j-1}$ for all $j\geq 1$. Therefore, for all $j\geq 0$, the support in $x$ of $a_j$ is contained in $U$. 
\end{proof}

Thus, for $s\in [-ch^\alpha,ch^\alpha]$ and $N \in \mathbb{N}$, we can construct the
$w_N^{ap}$ as in (\ref{WKBap}). We want to find $r_{h,N}$ such
that $w_N^{ap}$ satisfies (\ref{appeq}) and moreover to estimate
$r_{h,N}$ and $w_N^{ap}$. For this we start by estimating the
phase $\Phi$ and the amplitude $(a_j)_{j\in \mathbb{N}}$ as well
as their derivatives in $L^\infty$ norm.

\begin{prop}
\label{phf_estimates} For all $j,k \in \mathbb{N}$, $k\geq 1$ and $\beta \in
\mathbb{N}^d$ there exist constants $c_{k,\beta},
c_{k,\beta,j} >0$ such that functions $\Phi$ and $(a_j)_{j\in \mathbb{N}}$
constructed in Proposition \ref{time_exist} satisfy, for all $s\in
[-ch^\alpha, ch^\alpha]$, the estimates
\begin{equation}
\label{freq} \nX{\nabla_x^k \partial_\xi ^\beta
\Phi(s)}{L_x^\infty} \leq c_{k,\beta} h^{-\alpha\max(k-2,0)}
\end{equation}
and
\begin{equation}
\label{phase} \nX{\nabla_x^k \partial_\xi ^\beta
a_j(s)}{L_x^\infty} \leq c_{k,\beta,j} h^{-\alpha\max(k+j-1,0)}.
\end{equation}
Moreover, for $|\beta| \geq 2$,
\begin{equation}
\label{reffreq} 
\nX{\nabla_x \partial_\xi^\beta \Phi}{L_{x,\xi}^\infty} \leq ch^\alpha.
\end{equation}
\end{prop}

\begin{proof} In the proof of Proposition \ref{time_exist} we have deduced estimate (\ref{freq}) for $k=1$, $\beta=0$ : see (\ref{important}). Throughout this proof we consider $0<s\leq ch^\alpha$.

For $n \in \mathbb{N}$, $n\geq{0}$, we denote 
$$M_n(s)=\sup\limits_{|t|\leq s,\ |\gamma|=n} \nLx{\partial_x^\gamma \Phi(t)|}{\infty}.$$ 
Thus, estimate (\ref{important}) reads $M_1(s)\leq c$ for all $0<s\leq ch^\alpha$.

In order to estimate the functions $\nabla_x^k \partial_\xi^\beta \Phi$ for $k\geq 2$ or $|\beta|\geq 1$ we find the equations they verify by differentiating the equation (\ref{HJ}) satisfied by $\Phi$. We get that they satisfy transport equations. We estimate their $L^\infty$ norm by the transport Lemma \ref{trlemma} combined with induction on the order of derivatives. Having two parameters, we make first an induction on the order of derivatives in $x$, then in $\xi$. But first of all we need to estimate the $L^\infty$ norm of two derivatives in $x$ of $\Phi$, i.e. $M_2(s)$.

For $1\leq j,k \leq d$ let $V_{j,k}=\partial_{x_j}\partial_{x_k}\Phi$. Then $V_{j,k}$ verifies the equation
$$\partial_s V_{j,k} + f_l \partial_{x_l}V_{j,k} + F_{j,k} = 0, \ {V_{j,k}}(0)=0,$$
where we denote by $f_l=2G_h^{l,m}\partial_{x_m}\Phi$ and by $F_{j,k}$ the terms from $\partial_{x_j}\partial_{x_k} \left( G_h^{l,m}\partial_{x_l}\Phi \partial_{x_m}\Phi \right)$ except those that contain a 3-derivative in $\Phi$. We can decompose $F_{j,k}$ following the order of derivatives as follows $F_{j,k}=F_0+F_1+F_2$, where in $F_n$ there are $n$ derivatives on $G_h$. Combining estimates on $G_h$ (see Proposition \ref{regmetric}) with $M_1(s)\leq c$ we have :
$ \nX{F_0}{L_{x,s}^\infty} \leq cM_2(s)^2$, $\nX{F_1}{L_{x,s}^\infty}\leq cM_2(s)$ and $\nX{F_2}{L_{x,s}^\infty}\leq ch^{-\alpha}$. By the transport Lemma \ref{trlemma} we obtain $\nLx{V_{j,k}(s)} {\infty} \leq cs \nX{F_{j,k}}{L_{x,s}^\infty}$. Therefore,
$$M_2(s) \leq ch^{\alpha}(h^{-\alpha} + M_2(s) + M_2(s)^2).$$ We treat this inequation with a bootstrap method. Using that $M_2(0)=0$ we obtain $M_2(s) \leq c$ for all $s \leq ch^\alpha$.

Similarly, for $\gamma \in \mathbb{N}^d$ such that $|\gamma|\geq 3$, we denote by
$V_\gamma=\partial_x^\gamma \Phi$. By induction hypothesis $M_n(s) \leq ch^{-\alpha {\rm max}(n-2,0)}$, for all $n\leq |\gamma|-1$. Differentiating the Hamilton-Jacobi equation (\ref{HJ}) following $\partial_x^\gamma$, we get the transport equation verified by $V_\gamma$
$$\partial_s V_\gamma + f_l\partial_{x_l}V_\gamma + F_{\gamma} = 0.$$ 
Note that $f_l = 2 G_h^{l,m} \partial_{x_m}\Phi$ is the same for all $\gamma$'s and $F_{\gamma}$ equals $\partial_x^\gamma ( G_h^{l,m}\partial_{x_l}\Phi \partial_{x_m}\Phi)$ minus the terms that contain a $(|\gamma|+1)$-derivative in $\Phi$. Making a similar analysis with the one done for $M_2(s)$ we obtain $M_{|\gamma|}(s) \leq ch^\alpha M_{|\gamma|}(s) + ch^{-\alpha(|\gamma|-2)}$ and therefore $M_{|\gamma|}(s) \leq ch^{-\alpha(|\gamma|-2)}$.

In order to estimate the $L_x^\infty$ norm of $\partial_x^\gamma \partial_\xi^\beta \Phi$ we introduce $$M_{n,k}(s)=\sup_{|t|\leq s,\ |\gamma|=n,\ |\beta|=k} \nLx{\partial_x^\gamma \partial_\xi^\beta \Phi(t)}{\infty}.$$
Thus, estimate $M_n(s)\leq ch^{-\alpha \max(n-2,0)}$ reads $M_{n,0}(s)\leq ch^{-\alpha \max(n-2,0)}$ for all $n\geq 1$. We make a double induction : we increase $k$ by 1 and make a complete induction on $n \in \mathbb{N}$ as above. We obtain $M_{n,k}(s) \leq ch^\alpha M_{n,k}(s) + ch^{-\alpha(n-2)}$ and consequently $M_{n,k}(s) \leq ch^{-\alpha(|\gamma|-2)}$. 

Moreover, note that for $|\beta|\geq 2$ we have $\partial_{x_j}\partial_\xi^\beta \Phi(0)=0$ for all $1\leq j \leq d$ and therefore we obtain estimate (\ref{reffreq}), which reads $M_{1,|\beta|}\leq ch^\alpha.$

In a similar way we estimate the derivatives of $a_j$ in
$L^\infty$ norm. Note that for $j \geq 1$ the functions $a_j$ are
solutions of nonhomogeneous transport equations (\ref{trj}) with a source term that equals $i \tr_{G_h} a_{j-1}$. Thus, when we differentiate equation (\ref{trj}) with respect to
$x$ we get some powers of $h^{-\alpha}$ in the source term. This comes from the
frequency where we regularized the metric. And this loss explains why for bigger $j$'s we have a bigger loss in the $L^\infty$ norm of $a_j$.
\end{proof}

Let us recall that we denote by $v_0=\chi u_0$. For $N\in \mathbb{N}$ define
\begin{equation}
\label{hNrest} r_{h,N}=h^{N+2} \int_{\RR^d}
e^{i\frac{\Phi(s,x,\xi)}{h}} \tr_{G_h} a_N \widehat{v}_0 \left(
\frac{\xi}{h} \right) \frac{d\xi}{(2\pi h)^d}.
\end{equation}

Then $w_N^{ap}$ defined in (\ref{WKBap}) verifies, for $s\in
[-ch^\alpha, ch^{\alpha}]$ and $x\in \RR^d$, the equation

\begin{equation}
\left\{ \begin{array}{rcl}
i h\partial_s w_N^{ap} + h^2 \tr_{G_h} w_N^{ap} & = & r_{h,N}\\
{w_N^{ap}}_{|{s=0}} & = & F_h u_0.
\end{array}\right.
\end{equation}

\begin{prop}
\label{dispersive}
For $\alpha \geq \frac {1+r}{3+2r}$, where $r$ is an integer such that $r>\frac d2$, the approximate solution $w_N^{ap}$ constructed above satisfies, for
$s\in [-ch^\alpha, ch^\alpha]$, the following estimate
\begin{equation}
\label{appdisp} \nX{w_N^{ap}(s)}{L^\infty} \leq
\frac{c}{(h|s|)^{\frac d2}} \nL{v_0}{1}.
\end{equation}
\end{prop}

\begin{proof}
We write $w_N^{ap}$ as a kernel operator
$$w_N^{ap}(s,x)= \int_{\RR^d} K_h(s,x,y) v_0(y) dy,$$
where $K_h(s,x,y)=\int e^{i\frac{\Phi(s,x,\xi) - y\cdot \xi}{h}}
\sum_{j=0}^N h^j a_j(s,x,\xi) \frac{d\xi}{(2\pi h)^d}.$ Thus in order to control the absolute value of
$w_N^{ap}(s,x)$ by the norm $\nL{v_0}{1}$ it suffices to control the
$L^{\infty}_y$ norm of $K_h(s,x,y)$. The kernel is an oscillatory
integral, whose phase function can be written as
$$\Phi(s,x,\xi) - y\cdot \xi = (x-y)\cdot \xi + s \psi(s,x,\xi),$$
where $\psi$ is the remainder term from the Taylor expansion of $\Phi$ at first order
$$\psi(s,x,\xi)=\int_0^1 \partial_s \Phi(s\tau,x,\xi) {\rm d}\tau = -
\int_0^1 G_h^{l,m}(x) \partial_{x_l}\Phi(s\tau) \partial_{x_m}
\Phi(s\tau) {\rm d}\tau.$$ If we push the expansion to the second
order we get
\begin{equation} \label{etoile}
\psi(s,x,\xi) = -G_h^{l,m}(x)\xi_l \xi_m + s\int_0^1 (1-\tau)\partial_s^2 \Phi(s\tau) {\rm d}\tau.
\end{equation}
Setting $z=\frac{x-y}{s}$, we have $K_h(s)=\int e^{i\frac{s}{h}(z\cdot \xi +\psi(s,x,\xi))}
\sum_{j=0}^N h^j a_j(s,x,\xi) \frac{{\rm d}\xi}{(2\pi
h)^d}$ and we are interested in evaluating the $L^\infty_{z,x}$ norm of
$K_h$. Note that if $\frac{s}{h}$ is bounded we get immediately
that $|K_h(s)|\leq \frac{c}{h^d} \leq
\frac{c}{(s|h|)^\frac{d}{2}}$. Thus we can consider the rapport
$\lambda=\frac{s}{h}$ to be large. The kernel reads as
$$K_h(s)=\int e^{i\lambda F(s,z,x,\xi)}\sum_{j=0}^N h^j
a_j(s,x,\xi) \frac{{\rm d}\xi}{(2\pi h)^d},$$
with $F(s,z,x,\xi)= z\cdot \xi +\psi(s,x,\xi)$.

We want to apply the stationary phase lemma to estimate $K_h$. This
lemma says that the essential contribution in the integral must come
from points where the phase is stationary (critical nondegenerate
points). We shall use the stationary phase lemma under its simplest form
(lemma 7.7.3 in \cite{Ho}).
\begin{citelem}\textbf{\upshape{(\cite{Ho})}}
Let $A$ be a real symmetric non-degenerate matrix of dimension $d\times d$.
Then we have for every integer $k>0$ and integer $r > \frac{d}{2}$ :
\begin{equation}
\left| \int_{\RR^d} f(\xi)e^{i\lambda \frac{<A\xi,\xi>}{2}} {\rm d}\xi -\left(
\det \left(\frac{\lambda A}{2\pi i} \right)\right)^{-\frac{1}{2}} T_k(\lambda)
\right| \leq c_k \left( \frac{\norm{A^{-1}}{}}{\lambda} \right)^{\frac {d}{2}
+k} \sum_{|\beta| \leq 2k+r} \nLd{D^\beta f}
\end{equation}
for $f\in \mathcal{S}$ and for $T_k(\lambda)=
\sum_{j=0}^{k-1}(2i\lambda)^{-j} \frac{<A^{-1}D,D>^j f}{j!}(0)$.
\end{citelem}
For $s,x$ and $z$ fixed we want to show that the equation
$\partial_\xi F= 0$ has at most one solution. We write this equation
as $\xi = 2 G_h(x)^{-1}(z+s\int_0^1 (1-\tau)\partial_\xi \partial_s^2
\Phi(s\tau) {\rm d}\tau)$. It suffices to show that the right hand side is
contracting (as a function of $\xi$). For this we compute its
derivative with respect to $\xi$. Taking into account that
$s=ch^\alpha$ is small, it is enough to show that for $|\beta|=2$ the
expression $\left| 2G_h(x)^{-1}\int_0^1 (1-\tau )\partial_\xi^{\beta} \partial_s^2
\Phi(s\tau) {\rm d}\tau \right|$ is bounded (independently of $s,x,z,\xi$). We explicit $\partial_s^2 \Phi(s\tau)$ using the Hamilton-Jacobi equation (\ref{HJ}) verified by $\Phi$. The derivative with respect to $\xi$ distributes on the terms in
$\Phi$ and by using estimates (\ref{freq}) we conclude the boundedness of the expression. This implies the
uniqueness of the critical point of $F$. Let us call it $\xi_c$. The
phase $F$ decomposes as follows
$$F(\xi)=F(\xi_c) + \frac{1}{2}<\partial_\xi^2 F(\xi_c)(\xi -\xi_c),\xi-\xi_c> + R(\xi),$$
where the last term is the remainder term from the Taylor expansion at order 2
$$R(\xi)=3 \sum_{|\gamma|=3}\int_0^1 \partial_\xi ^\gamma
F(s,x,\xi_c+\theta(\xi-\xi_c))(1-\theta)^2 {\rm d}\theta \frac {(\xi -\xi_c)^\gamma}{\gamma!}.$$
We recall that $F(s,x,\xi) = z\cdot \xi -\int_0^1 G_h^{l,m}(x) \partial_{x_l} \Phi(s\theta) \partial_{x_m} \Phi(s\theta) {\rm d}\theta.$
Thus, for $|\gamma|=3$, in $\partial_\xi ^\gamma F$ at least two derivatives will bear on $\nabla_x\Phi(s\theta)$. By the refined estimates (\ref{reffreq}) we conclude that $|R(\xi)| \leq ch^\alpha |\xi - \xi_c|^3.$

We apply the stationary phase lemma for $k=1$ and $r>\frac{d}{2}$, $A_{i,j}=\partial_{\xi_i} \partial_{\xi_j}F(\xi_c)$  {\rm and} $f(\xi)=e^{i\lambda (F(\xi_c) + R(\xi))}a_0(s,x,\xi).$ We analyze the quantities that appear in the stationary phase lemma. Here $A$ is $O(h^\alpha)$ close to the regularized metric (see estimate (\ref{etoile})), so $\norm{A^{-1}}{}$ and $|{\rm det} A^{-1}|$ are bounded independently of $h$. We have $T_1(\lambda)=f(\xi_c)$ bounded, since $a_0$ is bounded. The only term we need to estimate in order to bound $|K_h(s,x)|$ is $\sum\limits_{|\beta|\leq 2+r}\nX{D_\xi^\beta f}{L^2_\xi}$. The function $f(\xi)=e^{i\lambda (F(\xi_c) + R(\xi))}a_0(s,x,\xi)$ being supported in $B(0,R)$, we have $\nL{D_\xi^\beta f}{2} \leq c\nL{D_\xi^\beta f}{\infty}$.
We explicit $D_\xi ^\beta f = \sum_{\eta \leq \beta} D_\xi^\eta \left( e^{i\lambda \left( F(\xi_c) + R(\xi) \right) }\right) D_\xi^{\beta - \eta}a_0.$
By a simple computation we get that $\nX{D_\xi^\eta R}{L_{x,\xi}^\infty} \leq ch^\alpha$, $\nX{D_\xi^\eta a_0}{L_{s,x,\xi}^\infty} \leq c$. Thus, the most important contribution in the sum comes from terms where the derivative bears on the exponential. Consequently,
$$\nL{D_\xi^\beta f}{\infty} \leq c \max (1, (\lambda h^\alpha) ^{|\beta|}).$$
From the stationary phase lemma we conclude that
\begin{equation}
\label{SPint2}
|K_h(s,x)| \leq \frac{c}{(\lambda h^2)^\frac d2} + c \sup_{|\beta|\leq 2+r} \frac{ \max (1, (\lambda h^\alpha) ^{|\beta|})}{h^d \lambda^{\frac d2 +1}}.
\end{equation}
As $\lambda h^2=|s|h$ and we are looking to prove
$|K_h(s,x)| \leq \frac{c}{(|s|h)^\frac d2},$
we want the second term of the sum to be small in front of the first one. This makes us to impose the following condition 
\begin{equation}
\label{strangecond}
\frac{\max (1,(\lambda h^\alpha)^{2+r})}{\lambda} \leq c,
\end{equation}
for all $0<h<1$. Let us recall that $\lambda =\frac{|s|}{h}$ is large and $|s|\leq ch^\alpha$.
Consequently, $ch\leq |s| \leq ch^\alpha$ and $\lambda h^\alpha=|s|h^{\alpha -1}$.

If $1>\alpha \geq \frac{1}{2}$ then $\lambda h^\alpha \leq ch^{2\alpha -1} \leq c$. Thus, $\max (1,(\lambda h^\alpha)^{2+r}) <c$ and (\ref{strangecond}) holds.

If $0<\alpha <\frac{1}{2}$ we have to study the case $ch \leq |s| \leq ch^{1-\alpha}$ and $ch^{1-\alpha} \leq |s| \leq ch^{\alpha}$. In the first case $\lambda h^\alpha \leq c$ and (\ref{strangecond}) holds as above. In the second case $\lambda h^\alpha \geq c$ and therefore condition (\ref{strangecond}) becomes $\frac{\max (1,(\lambda h^\alpha)^{2+r})}{\lambda} = |s|^{1+r}h^{-1-r}h^{\alpha(2+r)}<c$ for all $0<h<1$. Consequently $\alpha \geq \frac{1+r}{3+2r}$. 

We deduce from $\frac{1+r}{3+2r}\leq \alpha<1$ that (\ref{strangecond}) holds and combined with estimate (\ref{SPint2}) this implies $|K_h(s,x)| \leq \frac{c}{(|s|h)^\frac d2}$.

\end{proof}

\begin{prop}
\label{esthNrest}
The force term $r_{h,N}$ defined in (\ref{hNrest}) satisfies, for
$\sigma$ and $N \in \mathbb{N}$ such that $\sigma \leq N(1-\alpha)+2-d$ and for 
$s\in [-ch^\alpha, ch^\alpha]$, the estimate

$$\nX{r_{h,N}(s)}{H^\sigma(\RR^d)} \leq ch^{N(1-\alpha)+2-\alpha-\sigma-d}\nL{v_0}{1}.$$
\end{prop}

\begin{proof}
From (\ref{hNrest}) we deduce that
$$\nH{r_{h,N}}{\sigma} \leq h^{N+2} \nL{v_0}{1} \sup_y \int_{\RR^d} \nX{e^{i\frac{\Phi(s,x,\xi)-\xi y}{h}} \tr_{G_h} a_N(s,x,\xi)}{H_x^\sigma} \frac{{\rm d}\xi}{(2\pi h)^d} .$$
As $\nH{f}{\sigma} \leq \nL{f}{2} + \nL{D^\sigma f}{2}$ and the support in $x$ and $\xi$ of $a_N$ is compact (see Proposition \ref{time_exist}), we get
$$\nX{e^{i\frac{\Phi-\xi y}{h}} \tr_{G_h} a_N}{L_\xi^1(H_x^\sigma)} \leq c\nX{e^{i\frac{\Phi-y\xi}{h}} \tr_{G_h} a_N}{L_{x,\xi}^\infty} + c\nX{D^\sigma\left( e^{i\frac{\Phi-y\xi}{h}} \tr_{G_h} a_N \right) }{L_{x,\xi}^\infty}.$$
Note that when differentiating $e^{i\frac{\Phi-y\xi}{h}} \tr_{G_h} a_N$ once in $x$, the biggest contribution comes from differentiating the exponential and it is of order $h^{-1}$, while all the other terms contribute with at most $h^{-\alpha}$ growth. For $\sigma$ derivatives the order of magnitude is $h^{-\sigma}$. Thus $\nHx{e^{i\frac{\Phi-y\xi}{h}} \tr_{G_h} a_N}{\sigma} \leq ch^{-\alpha(N+1)-\sigma}$ and this uniformly in $s$ and $\xi$. Consequently,
$$\nHx{r_{h,N}}{\sigma} \leq ch^{N(1-\alpha)+2-\alpha-\sigma-d}\nL{v_0}{1}.$$
\end{proof}

\subsection{Strichartz inequality}

Further on, we consider $d=2,3$, as those are the only dimensions where we hope to get an existence theorem in the energy space from our Strichartz estimate. In this section we proceed to the proof of the Strichartz estimate as presented in section \ref{partition}. 

Let us recall the framework as introduced in section \ref{partition}. We have considered $(U_j,\kappa_j)_{j\in J}$ a finite covering with open charts of the manifold M. We have defined a family of spectral truncations on $M$: $J_h f = \sum_{j\in J}\tilde{\chi}_j \varphi(hD)(\chi_j(\kappa_j^{-1}) f(\kappa_j^{-1}))\kappa_j(x)$. We have generically denoted by $F_h$ a spectral truncation on functions of $\RR^d$, corresponding to one component of the partition of unity  $F_h f(y)=\tilde{\chi}(y) \varphi(hD) (\chi f)(y)$. In section \ref{ansatz} we have constructed $w_N^{ap}$ by the WKB method. This function verifies the Schrödinger equation for the regularized metric, with a small source term and initial data $F_h u_0$ (\ref{appeq}). From Proposition \ref{time_exist} we deduce $w_N^{ap}$ and $r_{h,N}$ are localized in the open chart corresponding to $\chi$. 

We resume the previous section in a lemma asserting that the function $w_N^{ap}$ constructed in (\ref{WKBap}) is close to $e^{it\tr_{G_h}} F_h u_0$ in $L_x^\infty$ norm. Moreover, $w_n^{ap}$ being localized in an open chart we can extend it to a function on the whole manifold. Thus, from the dispersive estimate on the approximate solution $w_N^{ap}$ in Proposition \ref{dispersive} we deduce a dispersive estimate for $e^{it\tr_{G_h}} J_h u_0$ on a small interval of time (of length $ch^{1+\alpha}$).

\begin{lem} \label{appsol} 
There exists a function $R_{h,N} :[-ch^{1+\alpha}, ch^{1+\alpha}]\times \RR^d \rightarrow \mathbb{C}$ such that
$$ e^{it\tr_{G_h}} F_h u_0 (x)=w_N^{ap}\left( \frac{t}{h},x\right)+ R_{h,N}(t,x)$$
and if we denote by $v_0 = \chi u_0$, then
$$\nL{R_{h,N}}{\infty}\leq ch^{N(1-\alpha) -1-d} \nL{v_0}{1}.$$
\end{lem}
\begin{proof}
The function $w_N^{ap}$ was constructed such that it satisfies the equation (\ref{appeq}). By the Duhamel formula applied to equation (\ref{appeq}), the following holds for all $s\in [-ch^\alpha, ch^\alpha]$: $$e^{ihs\tr_{G_h}} F_h u_0 = w_N^{ap} - i h^{-1}\int_0^s e^{ih(s-\tau) \tr_{G_h}}r_{h,N}(\tau)d\tau.$$
For $t\in [-ch^{1+\alpha},ch^{1+\alpha}]$, we denote by $$R_{h,N}(t,x) = -i h^{-2} \int_0^t e^{i(t-l)\tr_{G_h}} r_{h,N}\left( \frac {l}{h}, x\right) {\rm d}l.$$ Using the change of variable $h\tau=l$ we conclude the first identity holds. 

Let us estimate $\nL{R_{h,N}(t)}{\infty}\leq h^{-2}\int_0^t \nL{e^{i(t-l)\tr_{G_h}} r_{h,N} (\frac{l}{h})} {\infty}\!\! {\rm d}l$. Using the Sobolev imbedding $H^{2} \subset L^\infty$, as $2>\frac{d}{2}$, we have to estimate the $H_x^2$ norm of $e^{i(t-l)\tr_{G_h}} r_{h,N}(\frac{l}{h})$. For this we need to commute $\tr$ with $e^{it\tr_{G_h}}$. We use the following elliptic regularity lemma (see e.g. \cite{GbgTr})

\begin{citelem}For all $u\in L^2(\RR^d)$ such that $\tr_{G_h}\in L^2(\RR^d)$ we know that $u\in H^{2}(\RR^d)$ and the following estimate holds
\begin{equation} \label{paire}
\nH{u}{2}\leq c(\nL{u}{2}+ \nL{\tr_{G_h} u}{2}) \leq c\nH{u}{2}.
\end{equation}
\end{citelem}
\noindent Consequently, $$\nLx{R_{h,N}}{\infty} \leq ch^{-2}\int_0^t \nLx{e^{i(t-l)\tr_{G_h}} r_{h,N} \left( \frac{l}{h} \right)}{2} + \nLx{e^{i(t-l)\tr_{G_h}} \tr_{G_h}r_{h,N} \left(\frac{l}{h}\right)}{2} {\rm d}l.$$ Using the conservation of the $L^2$ norm by the flow $e^{it\tr_{G_h}}$ and the second inequality from (\ref{paire}), we obtain
$$\nLx{R_{h,N}(t)}{\infty} \leq ch^{-2} \int_0^t \nHx{r_{h,N}\left(\frac lh \right)}{2}{\rm d}l.$$
Thus, using the estimate of the remainder term $r_{h,N}$ seen in Proposition \ref{esthNrest} for $\sigma=2$ and $|t|\leq ch^{1+\alpha}$, the result follows.
\end{proof}

\begin{prop}For all $u_0 \in L^1(M)$ there exists constants $C>0$ and $c>0$ such that, for all $t\in [-ch^{1+\alpha}, ch^{1+\alpha}]$, the following dispersive estimate holds
\begin{equation}
\label{hregdispersion}
\nLx{e^{it\tr_{G_h}} J_h u_0}{\infty} \leq \frac{C}{|t|^\frac d2} \nL{u_0}{1}.
\end{equation}
\end{prop}

\begin{proof} Let us recall that in local coordinates $J_h$ is a sum of truncations $F_{h,j}$ corresponding to $\chi_j$. From Lemma \ref{appsol} and the semiclassical dispersive estimate (\ref{appdisp}) we obtain
$$\nLx{e^{it\tr_{G_h}}F_{h,j} u_0}{\infty} \leq \frac{c}{|t|^{\frac{d}{2}}} \nL{v_{0,j}}{1} + ch^{N(1-\alpha) -1-d} \nL{v_{0,j}}{1}.$$ As $\sum_{j\in J}\chi_j=1$ and $0\leq \chi_j \leq 1$, we can sum both left and right side terms. We obtain
$$\nLx{e^{it\tr_{G_h}}J_h u_0}{\infty} \leq \frac{c}{|t|^{\frac{d}{2}}} \nL{u_0}{1} + ch^{N(1-\alpha) -1-d} \nL{u_0}{1}.$$
For $0<\alpha<1$, we can find $N\in \mathbb{N}$ such that $N(1-\alpha)-1-d >0$. Therefore, the second term is absorbed by the first one and the result follows.

\end{proof}

Having a dispersive estimate $L^1\rightarrow L^\infty$ we obtain the following spectrally truncated Strichartz estimate (as well as its adjoint form).
\begin{prop} For all couples $(p,q)$ which are admissible in dimension $d$ and $I_h$ an interval of time such that $|I_h|=ch^{1+\alpha}$, we have
\begin{equation}
\label{hStrinter}
\nX{J_h^* e^{it\tr_{G_h}} u_0}{L^p(I_h,L^q(M))} \leq c\nL{u_0}{2}
\end{equation}
and
\begin{equation}
\label{hStrinter*}
\nL{\int_{I_h}e^{it\tr_{G_h}}J_h F(t,x){\rm d}t}{2} \leq c\nX{F}{L^p(I_h,L^q(M))}.
\end{equation}
\end{prop}

\begin{proof}
This is quite a straightforward result from the following $TT^*$ method (which was optimized by Keel and Tao \cite{KT} for the endpoint case).

\begin{citelem} A parametrized family of operators $U(t) : L^2 \rightarrow L^2$ that obeys, for all $t$, the energy estimate
$$\nLx{U(t)f}{2} \leq c \nLx{f}{2}$$
and the decay estimate
$$\nLx{U(t) U^*(s)f}{\infty} \leq \frac{c}{|t-s|^{\frac d2}} \nLx{f}{1}$$
satisfies, for all admissible pairs $(p,q)$, $(p_1,q_1)$ in dimension $d$, the estimates

$$\nLtLx{U(t)f}{p}{q} \leq B_1(q)\nL{f}{2},$$
$$\nL{\int U^*(s)F(s)ds}{2} \leq B_2(q) \nLtLx{F}{\bar{p}}{\bar{q}}.$$
\end{citelem}
\noindent We consider the operator $U_h(t)=J_h^* e^{it\tr_{G_h}}$. 
Thus $U_h(t)U_h^*(s)u_0 = J_h^* e^{i(t-s)\tr_{G_h}}J_h u_0$.
We use the boundedness of $J_h$ on $L^p$ spaces for all $1\leq p \leq \infty$ (see Lemma \ref{commJhGh}) to conclude from inequality (\ref{hregdispersion}) that $U_h$ satisfies the decay estimate as requested by the $TT^*$ method.
\end{proof}

\begin{rem} Let us suppose that in estimate (\ref{hStrinter}) we have $e^{it\tr_G}$ instead of $e^{it\tr_{G_h}}$. Still, we could not sum over all frequencies as on the left side there is a term that does not depend on the frequency.
\end{rem}
In the following we deduce an Strichartz inequality that will sum on all frequencies. Let $\varphi$ be as in (\ref{LPdecomp}) and $\tilde{\varphi}$ supported on an annulus such that $\tilde{\varphi} \equiv 1$ on an open neighborhood of the range of $\nabla_x \Phi$ near the support of $a_0$. From (\ref{important}) we know $c\leq |\nabla_x \Phi(s,x,\xi)|\leq C$. Moreover, as $\nabla_x\Phi(0,x,\xi)=\xi$ for $\xi\in supp(\varphi)$, we conclude that $\tilde{\varphi}\equiv 1$ on the support of $\varphi$. Let $\tilde{J}_h^* = \sum \chi_j(x) \tilde{\varphi}(hD) \tilde{\chi}_j$.

\begin{prop}
\label{intermediar} For all $u_0\in L^1(M)$ there exists $R_N:[-ch^{1+\alpha}, ch^{1+\alpha}]\times \RR^d \rightarrow \mathbb{C}$ such that for all $|t|\leq ch^{1+\alpha}$
\begin{equation}
\label{egorov}
\tilde{J}_h^* e^{it\tr_{G_h}} J_h u_0 (x) = e^{it\tr_{G_h}} J_h u_0(x) + R_{N}(t,x),
\end{equation}
and for all $N_0>0$ we can choose $N$ such that , $\nL{R_{N}(t)}{\infty}\leq h^{N_0}\nL{u_0}{1}$.
\end{prop}
Note that this proposition states the localization of the flow at the same frequency as the initial data on a time scale $h^{1+\alpha}$.

\begin{proof}
Let us recall that $\tilde{J}_h^* e^{it\tr_{G_h}} J_h u_0  = \sum_{j,l\in J}\tilde{F}_{j,h}^* e^{it\tr_{G_h}} F_{l,h} u_0.$
We pass into semiclassical coordinates by setting $t=hs$ and use the WKB approximation (as resumed by Lemma \ref{appsol}) to express $$ \tilde{F}_{j,h}^* e^{ihs\tr_{G_h}} F_{l,h} u_0 (x)= \chi_j (x) \frac{1}{h^d}\int_{\RR^d} \tilde{\rho} \left(\frac{x-y}{h} \right) \tilde{\chi}_j(y) \left( w_N^{ap}(s,y) + R_{h,N}(hs,y) \right){\rm d}y.
$$
We make the change of variable $y=x-hz$. We denote by $a_{h,N}(s,x,\xi)=\sum_{k=1}^N a_k(s,x,\xi)$ the amplitude of the WKB ansatz. We make a Taylor expansion in $x$ following $hz$. Thus, the main part of $\tilde{F}_{j,h}^* e^{ihs\tr_{G_h}} F_{l,h} u_0 (x)$ reads
$$
\chi_j(x) \int_{\RR^d \times \RR^d} \tilde{\rho}(z) \tilde{\chi}_j(x) e^{\frac ih \Phi(s,x,\xi)} e^{-iz\cdot \nabla_x \Phi(s,x,\xi)} a_{h,N}(s,x,\xi) \hat{v}_0\left(\frac{\xi}{h}\right) {\rm d}z \frac{{\rm d}\xi}{(2\pi h)^d}
$$
Using that $\int_{\RR^d } \tilde{\rho}(z) e^{-iz\cdot \nabla_x \Phi(s,x,\xi)} {\rm d}z = \tilde{\varphi}(\nabla_x \Phi) (s,x,\xi)$ and the hypothesis $\tilde{\varphi}(\nabla_x \Phi)\equiv 1$, we obtain that it equals $\chi_j (x) w_N^{ap}(s,x)$. We apply again Lemma \ref{appsol} and get that $$\tilde{F}_{j,h}^* e^{ihs\tr_{G_h}} F_{l,h} u_0 (x) = \chi_j (x) e^{ihs\tr_{G_h}} F_{l,h} u_0 (x) + \chi_j (x) R_{h,N}(hs,x).$$

Moreover, for $|\beta|\geq 1$ we have $\int_{\RR^d } z^\beta \tilde{\rho}(z) e^{-iz\cdot \nabla_x \Phi(s,x,\xi)} {\rm d}z = (\partial^\beta \tilde{\varphi} )(\nabla_x \Phi) (s,x,\xi)=0.$ Thus, we get that all the terms from the Taylor and WKB expansion are null, except those containing some remainder terms. 

We denote by $R_N$ the sum of the remainders from WKB approximation and Taylor expansions. Both kinds of remainders contain a sufficiently large power of $h$ to be treated as in Proposition \ref{esthNrest} and Lemma \ref{appsol}.
\end{proof}

From Proposition \ref{intermediar} we can easily deduce a Strichartz estimate similar to (\ref{hStrinter}) that would sum over all frequencies. Nevertheless, estimating the difference $(e^{it\tr_{G_h}} - e^{it\tr_G})J_h u_0$ in the $L^p(I_h,L^q(M))$ norm turns out to be a difficult task, as we know very little on $e^{it\tr_G}$. We prefer to estimate the $L^p(I_h,L^q(M))$ norm of $J_h^* (e^{it\tr_{G_h}} - e^{it\tr_G}) u_0$. In view of this, we deduce from Proposition \ref{intermediar}, the following result.

\begin{prop} For $u_0\in H^1(M)$ and $(p,q)$ an admissible pair in dimension $d=2$ or $3$, the following holds, for $|I_h|=ch^{1+\alpha}$,
 
\begin{equation}
\label{hStr}
\nLiLM{J_h^* e^{it\tr_{G_h}} u_0}{p}{I_h}{q} \leq c h\nH{u_0}{1}.
\end{equation}
\end{prop}

\begin{proof}
We use an adjoint argumentation. Let $F\in L^{p'}(I_h,L^{q'})$. Then
$$<J_h^* e^{it\tr_{G_h}} u_0, F(t,x)>_{L_t^p(L_x^q),L_t^{p'}(L_x^{q'})} = <u_0,\int_{I_h} e^{-it\tr_{G_h}}J_h F(t,x){\rm d}t>_{L_x^2}.$$
We apply (\ref{egorov}) for $u_0=F(t,\cdot)$ and thus the previous expression equals
$$<\tilde{J}_h u_0,\int_{I_h} e^{-it\tr_{G_h}}J_h F(t,x){\rm d}t>_{L_x^2} + <u_0,\int_{I_h} R_N( F(t,x)){\rm d}t>_{L_x^2}.$$
Using the Strichartz inequality under its adjoint form (\ref{hStrinter*}) and the estimates on $R_N$ from (\ref{egorov}), we obtain
$$\left| <J_h^* e^{it\tr_{G_h}} u_0, F(t,x)>_{L_t^p(L_x^q),L_t^{p'}(L_x^{q'})} \right| \leq \nL{\tilde{J}_h u_0}{2} \nLL{F}{p'}{q'} + h^N \nL{u_0}{2} \nLL{F}{p'}{q'}$$
and the result follows from $\nL{\tilde{J}_h u_0}{2}\leq ch\nH{u_0}{1}$.
\end{proof}

Note that estimate (\ref{hStr}) sums for $h=2^{-k}$, $k\in \mathbb{N}$. As we are looking for a Strichartz inequality for $e^{it\tr_G}$, before summing, we will estimate the $L^p(I_h,L^q)$ norm of the difference
\begin{equation}
\label{Fhrest}
R(t)u_0 = J_h^* e^{it\tr_{G_h}}u_0 - J_h^* e^{it\tr_G}u_0.
\end{equation}

We have already introduced the notation $J_h^* = \tilde{J}_h^* J_h^* +T_h$ and estimated $\nX{T_h}{L^p \rightarrow L^p} \leq c_N h^N$ in Lemma \ref{Rlemma}. We use it here in order to write 
\begin{equation}
\label{Rdecomp}
R(t)=\tilde{J}_h^* R(t) + T_h (e^{it\tr_{G_h}}u_0 - e^{it\tr_G}u_0).
\end{equation}

\begin{prop}
The operator $R(t)u_0$ defined in (\ref{Fhrest}) satisfies, for all admissible pairs $(p,q)$ and $|I_h|=ch^{1+\alpha}$, to
\begin{equation}
\label{estFhrest} \nLiLM{R(t)u_0}{p}{I_h}{q} \leq c h^{2\alpha} \nH{u_0}{1}.
\end{equation}
\end{prop}

\begin{proof}
We bound the last term from (\ref{Rdecomp}) using estimate (\ref{RRest}) and the Sobolev imbedding $H^1\subset L^q$, for $2\leq q \leq \frac{2d}{d-2}$ (strict inequality for $d=2$) :
$$\nLx{T_h (e^{it\tr_{G_h}}u_0 - e^{it\tr_G}u_0)}{q} \leq c h^2 \nLx{e^{it\tr_{G_h}}u_0 - e^{it\tr_G}u_0}{q}\leq ch^2 \nH{u_0}{1}.$$
Consequently, $$\nLiLM{T_h (e^{it\tr_{G_h}}u_0 - e^{it\tr_G}u_0)}{p}{I_h}{q} \leq ch^{2+\frac{1+\alpha}{p}}\nH{u_0}{1}.$$ 
By a simple computation one can see that $R(t)u_0$ verifies the following equation
\begin{equation}
\left\{
\begin{array}{rcl}
(i \partial_t  + \tr_{G_h}) R(t)u_0 & = & [\tr_{G_h},J_h^*]\left( e^{it\tr_{G_h}}-e^{it\tr_G}\right)u_0 + J_h^*(\tr_G-\tr_{G_h})e^{it\tr_G}u_0,\\
R(0)u_0 & = & 0.
\end{array}\right.
\end{equation}
By the Duhamel formula we get that $R(t)u_0$ equals
$$
\int_0^t e^{i(t-\tau)\tr_{G_h}}[\tr_{G_h},J_h^*]\left( e^{i\tau \tr_{G_h}}-e^{i\tau \tr_G}\right)u_0 {\rm d}\tau + \int_0^t e^{i(t-\tau)\tr_{G_h}}J_h^*(\tr_G-\tr_{G_h})e^{i\tau \tr_G}u_0 {\rm d}\tau.
$$
We decompose $\tilde{J}_h^* R(t)=I_1+I_2,$
where $$I_1=\int_0^t \tilde{J}_h^* e^{i(t-\tau)\tr_{G_h}}[\tr_{G_h},J_h^*]\left( e^{i\tau \tr_{G_h}}-e^{i\tau \tr_G}\right)u_0 {\rm d}\tau$$ and $$I_2=\int_0^t \tilde{J}_h^* e^{i(t-\tau)\tr_{G_h}}J_h^*(\tr_G-\tr_{G_h})e^{i\tau \tr_G}u_0 {\rm d}\tau.$$
We apply the Minkowski inequality (as $p \geq 2$) as follows 
$$
\begin{array}{rcl}
\nLiLM{I_1}{p}{I_h}{q}&\leq &\int_0^T \nLiLM{1_{\tau\leq t} \tilde{J}_h^* e^{i(t-\tau)\tr_{G_h}}[\tr_{G_h},J_h^*]\left( e^{i\tau \tr_{G_h}}-e^{i\tau \tr_G}\right)u_0}{p}{I_h}{q}\!\!\!\!\!\!\!{\rm d}\tau \\
&\leq & \int_0^T \nL{[\tr_{G_h},J_h^*]\left( e^{i\tau \tr_{G_h}}-e^{i\tau \tr_G}\right)u_0}{2}{\rm d}\tau \leq ch^{1+\alpha} \nH{u_0}{1},
\end{array}
$$ 
where we have used the Strichartz estimate (\ref{hStrinter}), Lemma \ref{commJhGh} and the $H^1$ conservation law of both $e^{it\tr_G}$ and $e^{it\tr_{G_h}}$. 

Similarly, we estimate $\nLiLM{I_2}{p}{I_h}{q}\leq ch^{2\alpha}\nH{u_0}{1}$ using Lemma \ref{diffGGh}. 

For $\frac{1+r}{3+2r} \leq \alpha <1$ (see Proposition \ref{dispersive}), the minimum of $2+\frac{1+\alpha}{p}$, $1+\alpha$ and $2\alpha$ is $2\alpha$ and the result follows.
\end{proof}

We are now ready to deduce a Strichartz inequality of the spectrally truncated flow on a small time interval.

\begin{prop}
\label{PFIhStr}
For $\varphi$ a $C^\infty$ function supported in an annulus and $u_0 \in H^1$, for each admissible pair $(p,q)$ in dimension $d=2$ or $3$ and for each interval of time $I_h$, $|I_h|=h^{1+\alpha}$, the following Strichartz inequality holds
\begin{equation}
\label{FIhStr}
\nLiL{J_h^*e^{it\tr_G}u_0}{p}{I_h}{q} \leq c h^{\min (1,2\alpha)} \nH{u_0}{1}.
\end{equation}
Moreover, for $M$ flat outside a compact set, estimate (\ref{FIhStr}) also holds for $J_{h,\infty}=J_h + F_{1,\infty} + F_{2,\infty}$.
\end{prop}

\begin{proof}
Combining the estimate on the remainder term (\ref{estFhrest}) with the Strichartz inequality of the spectrally truncated flow for the regularized metric (\ref{hStr}) we get (\ref{FIhStr}).

For the case $M$ exterior of a compact set, let $v(t,x)=F_{\infty} e^{it\tr_G}u_0$, where $F_\infty$ defined by (\ref{LPdecompFlat}). Then $v$ satisfies
$$
\left \{
\begin{array}{rcl}
i \partial_t v + \triangle v &=& [\tr, F_{\infty}]e^{it\tr_G}u_0\\ 
v_{|_{t=0}} &=& F_{\infty} u_0,
\end{array}
\right.
$$
the Schrödinger equation with standard Laplacian on $\RR^d$. Therefore, we can apply the classical Strichartz inequality to the Duhamel formula 
$$v(t,x)=e^{it\tr}F_{\infty} u_0 + \int_0^t e^{i(t-s)\tr}[\tr,F_{\infty}] e^{is \tr_G}u_0 {\rm d}\tau.$$ Note that $[\tr,F_{\infty}]$ is a bounded $H^1$ to $L^2$ operator. Thus, we obtain

\begin{equation}
\label{Strinfty}
\nLiL{F_{\infty} e^{it\tr_G}u_0}{p}{I_h}{q} \leq c h\nH{u_0}{1}.
\end{equation}
Note that this estimate is also true on an interval of time of length $ch$, but the estimate on $I_h$, $|I_h|=ch^{1+\alpha}$, is all we need. 
\end{proof}

\begin{rem} In the following we shall use Proposition \ref{PFIhStr} to obtain the Strichartz inequality for $e^{it\tr_G}u_0$ on $M$. For $M$ flat outside a compact set one needs to replace $J_h$ by $J_{h,\infty}$.
\end{rem}

We want to have similar results on a fixed time interval. Knowing the conservation of the $H^1$ norm by the flow $e^{it\tr_G}$, one can sum the results on small intervals adjacent to each other.

\begin{prop}
For $\varphi$ a $C^\infty$ function supported in an annulus and $u_0 \in H^1$, for $(p,q)$ an admissible pair, the following inequality holds
\begin{equation}
\label{FStr}
\nLiL{J_h^*e^{it\tr_G}u_0}{p}{[0,1]}{q} \leq c h^{\gamma- \frac{1 + \alpha}{p}} \nH{u_0}{1},
\end{equation}
where $\gamma=\min (1,2\alpha)$. 
\end{prop}

\begin{proof}
We write the interval $[0,1]$ as an union of intervals $$[0,1] = \cup _{l \in L} I_h^l,$$ where $I_h^l=[t_l,t_{l+1}]$, $0 \leq t_{l+1}-t_l \leq c h^{1+\alpha}$ and $\#{L}=ch^{-1-\alpha}$. Thus, on each interval $I_h^l$, inequality (\ref{FIhStr}) holds
$$\nLiL{J_h^*e^{it\tr_G}u_0}{p}{I_h^l}{q} \leq c h^\gamma \nH{e^{i t_l \tr_G}u_0}{1}.$$
We can sum the $p$th power of those inequalities. Using the conservation of the $H^1$ norm by the flow $e^{it\tr_G}$, we get
$\sum_{l \in L}\nH{e^{i t_l \tr_G}u_0}{1}^p \leq c h^{-1 -\alpha} \nH{u_0}{1}^p$. Consequently,
$$\nLiL{J_h^*e^{it\tr_G}u_0}{p}{[0,1]}{q} = \left( \sum_{l\in L} \nLiL{J_h^*e^{it\tr_G}u_0}{p}{I_h^l}{q}^p \right) ^{\frac 1p} \leq ch^{\gamma -\frac{1+\alpha}{p}} \nH{u_0}{1}.$$
\end{proof}
 
Having a Strichartz inequality for the spectrally truncated flow $J_h^* e^{it\tr_G}u_0$ on a fixed time interval, we take the sum for $h=2^{-k}$ for $k\in \mathbb{N}$ of those inequalities. Let $\varphi_0$ and $\varphi$ be like in (\ref{LPdecomp}). 

\begin{rem} One way of summing is to apply the triangle inequality to (\ref{2LPdecomp}) in order to estimate the $L^p(I,L^q)$ norm of the flow $e^{it\tr_G}u_0$ using the estimate (\ref{FStr})
$$
\begin{array}{rcl}
\nLiLM{e^{it\tr_G}u_0}{p}{I}{q} &\leq & \nLiLM{J_0 e^{it\tr_G}u_0}{p}{I}{q} + \sum\limits_{k \in \mathbb{N}} \nLiLM{J_{2^{-k}}e^{it\tr_G}u_0}{p}{I}{q}\\
& \leq & c\nH{u_0}{1} + \sum\limits_{k \in \mathbb{N}} 2^{-k(\gamma - \frac{1+\alpha}{p})} \nH{u_0}{1} \leq c \nH{u_0}{1}.
\end{array}
$$

For summing the terms in the right hand side we used $\gamma - \frac{1+\alpha}{p} > 0$, which is always true. Doing so we did not gain with respect to Sobolev imbeddings. In fact, using the admissibility condition $\frac 2p + \frac dq = \frac d2$, the Sobolev imbedding $H^{\frac 2p}(M) \hookrightarrow L^q(M)$ holds and we can trivially obtain
$$\nLiLM{e^{it\tr_G}}{p}{I}{q} \leq c \nX{e^{it\tr_G}u_0}{L^\infty(I,H^{\frac 2p})} \leq c \nX{e^{it\tr_G}u_0}{L^\infty(I,H^{1})} \leq c \nH{u_0}{1}.$$
\end{rem}

Having an $H^1$ norm of $u_0$ in the right hand side term we try to improve the way of summing on the left hand side. We denote by $W^{\sigma,q}(M)$ the domain of $(1-\tr)^{-\frac{\sigma}{2}}$ in $L^q(M)$
$$W^{\sigma,q} = \{f \in L^q\ s.t.\ (1-\tr)^{\frac{\sigma}{2}} f\in L^q \}$$ endowed with the norm
$$\nX{f}{W^{\sigma,q}}= \nL{(1-\tr)^{\frac{\sigma}{2}}f}{q}.$$
Let $\epsilon >0$ be a small parameter such that if we denote by $$\sigma(\alpha)=\gamma - \frac{1+\alpha}{p} - \epsilon$$ we have $\sigma (\alpha)>0$. Note that $0<\alpha<1$ and $\gamma=\min (1,2\alpha)$ imply $0<\sigma(\alpha)<1$.

We bound the $L^p(I,W^{\sigma(\alpha),q})$ norm of the flow $e^{it\tr_G}$ by the $H^1$ norm of the initial data for all $\frac{1+r}{3+2r}\leq \alpha \leq 1$ and then prove that the best estimate is obtained for $\alpha=\frac 12$.

\begin{prop}
Let $I$ be a finite time interval and $(p,q)$ an admissible pair in dimension $d$. Then for all $\epsilon>0$ small, there exists a constant $c>0$ such that for all $u_0 \in H^1(M)$ the following holds
\begin{equation}
\label{SSStrest}
\nX{e^{it\tr_G}u_0}{L^p(I,W^{1-\frac {3}{2p}-\epsilon,q}(M))} \leq c \nHM{u_0}{1}.
\end{equation}
\end{prop}

\begin{proof}
Let $\epsilon>0$ such that $\sigma(\alpha)=\gamma - \frac{1+\alpha}{p}-\epsilon>0$.
As above, using (\ref{LPdecomp}), we have $$\nX{f}{W^{\sigma,q}} \leq \nL{J_0 f}{q}+ \sum_{k \in \mathbb{N}} 2^{k\sigma} \nL{J_{2^{-k}}f}{q}.$$ From estimates (\ref{FStr}), we obtain
$$
\begin{array}{lcl}
\nX{e^{it\tr_G}u_0}{L^p(I,W^{\sigma,q})} &\leq & \nLiLM{J_0 e^{it\tr_G}u_0}{p}{I}{q} + \sum_{k \in \mathbb{N}} 2^{k\sigma}\nLiLM{J_{2^{-k}}e^{it\tr_G}u_0}{p}{I}{q}\\
& \leq & c\nH{u_0}{1} + \sum\limits_{k \in \mathbb{N}} 2^{-k \epsilon} \nH{u_0}{1} \leq c \nH{u_0}{1}.
\end{array}
$$

Having a fixed norm $\nH{u_0}{1}$ on the right hand side, we want to find the best norm on the left hand side. As we have seen, we have the estimate with $\sigma=0$ for free. Thus, we want to find the largest $\sigma>0$ that satisfies. We analyze the function $\sigma(\alpha)=\min(1,2\alpha)-\frac{1+\alpha}{p}-\epsilon$ for $\alpha \in [\frac37, 1)$. Let us recall that the inferior bound comes from Proposition \ref{dispersive} $\left( \alpha\geq \frac{1+r}{3+2r} \right)$ applied to $r=2>\frac d2$. As $\sigma(\alpha)$ increases for $\alpha\leq \frac 12$ and decreases for $\alpha \leq \frac12$, we obtain that the function $\sigma(\alpha)$ takes its maximal value for $\alpha=\frac 12$ and it equals $1-\frac {3}{2p}-\epsilon$.
\end{proof}

We are now ready to deduce the result of Theorem \ref{StrThm}.
\begin{proof} \emph{of Theorem \ref{StrThm}}
From the elliptic regularity of $\tr_{G}$ (as in estimate (\ref{paire})) we know $\nL{(1-\tr)^\frac 12 u_0} {q}\approx \nL{(1-\tr_{G})^\frac 12 u_0}{q}$. Using the complex interpolation method we obtain it for fractional powers $\frac{\sigma}{2}$ for $0<\sigma<1$. Consequently,
$$\nL{(1-\tr)^{\frac 12 (1-\frac{3}{2p}-\epsilon)} u_0} {q}\approx \nL{(1-\tr_{G})^{\frac 12 (1-\frac{3}{2p}-\epsilon)} u_0}{q}.$$ This can also be read as follows: for all $u_0\in W^{1-\frac{3}{2p}-\epsilon,q}$ there exists $f\in L^q$ such that $u_0=(1-\tr_{G})^{-\frac 12 (1-\frac{3}{2p}-\epsilon)}f$. We introduce it into estimate (\ref{SSStrest}) and using that $e^{it\tr_G}$ commutes with $(1-\tr_{G})^{-\frac 12 (1-\frac{3}{2p}-\epsilon)}$, we obtain
$$\nLiLM{e^{it\tr_G}f}{p}{I}{q}\leq c\nH{f}{\frac{3}{2p}+\epsilon}.$$
\end{proof}

In order to control the nonlinear term in the proof of the local existence (see Theorem \ref{LEThm}) we have assumed and used the $L^p(L^\infty)$ norm estimate of the flow (corollary \ref{LStr}). We deduce it from estimate (\ref{SSStrest}) using the Sobolev imbeddings. Those state that for $\sigma$, $q$ and $d$ such that $\sigma q > d$ we have $W^{\sigma,q}(M) \hookrightarrow L^\infty(M)$. We want to find in which dimension we can deduce the control of the $L^p(L^\infty)$ norm. We combine the admissibility condition $\frac dq = \frac d2 - \frac 2p$, $(p,q,d) \neq (2,\infty,2)$ with the Sobolev condition for $\sigma = 1 -\frac{3}{2p}-\epsilon$. This yields the condition
\begin{equation}
\label{Sobolevcond}
1-\frac {3}{2p}- \epsilon > \frac d2 - \frac 2p.
\end{equation}
Consequently, $d\leq 2$ and this proves corollary \ref{LStr}. For $d=3$ the Strichartz inequality (\ref{Strest}) does not give us control of the $L^p(L^\infty)$ norm.

To  our knowledge, in the case of domains of $\mathbb{R}^3$, a local existence result in $H^1$, for instance for a cubic nonlinearity ($\beta=2$), remains an open problem.

\vskip 2mm

\noindent {\bf Acknowledgments :} \textit{The author would like to thank D.Tataru who gave the initial idea of this work. She would also like to thank P.Gérard for guidance from idea to achievement. This result is part of author's PhD thesis in preparation at Université Paris Sud, Orsay, under P.Gérard's direction.}

\textit{This work has started during the stay of the author at UC Berkeley in fall 2002, supported by a grant from France Berkeley Funds.}

\end{document}